\documentclass[12pt]{iopart}
\usepackage{amsthm}
\usepackage{cite}  

\bibstyle{unsrt}
\usepackage[colorlinks,linkcolor=black,anchorcolor=black,citecolor=black]{hyperref}
\usepackage{amssymb}
\usepackage{graphicx} 
\usepackage{float} 
\usepackage{subfigure}
\usepackage{xcolor}
\eqnobysec

\newtheorem{theorem}{Theorem}[section]

\newcommand{\mz}{\mathbf{z}}
\newcommand{\my}{\mathbf{y}}
\newcommand{\mx}{\mathbf{x}}
\newcommand{\mc}{\mathbf{c}}
\newcommand{\mb}{\mathbf{b}}
\newcommand{\mf}{\mathbf{f}}

 		\newcommand{\R}{\mathbb{R}}		
 		
 \newcommand{\Rn}{\R^n} 		   		\newcommand{\Rnn}{\R^{n\times n}}   		\newcommand{\Rnm}{\R^{n\times m}}  
 \newcommand{\Rmn}{\R^{m\times n}}  
 \newcommand{\Rmm}{\R^{m\times m}}

 \newcommand{\ct}{\centering}

 \newcommand{\bmt}{\left[\begin{matrix}}    
 \newcommand{\emt}{\end{matrix}\right]}

 \newcommand{\aq}{\alpha_{1}}  
 \newcommand{\aw}{\alpha_{2}}  
  \newcommand{\bq}{\beta_{1}}  
 \newcommand{\bw}{\beta_{2}}

 \newcommand{\Mbq}{\mathcal{M}_{\beta_1}}

\begin{document}

\title[Parameter identification from single trajectory data]{Parameter identification from single trajectory data:  from linear to nonlinear}
\author{X Duan$^1$, J E Rubin$^1$, D Swigon$^1$}
\address{$^1$ Department of Mathematics, University of Pittsburgh, Pittsburgh, PA 15260}
\ead{xid33@pitt.edu, jonrubin@pitt.edu, swigon@pitt.edu}

\begin{abstract}
Our recent work lays out a general framework for inferring information about the parameters and associated dynamics of a differential equation model from a discrete set of data points collected from the system being modeled.  Rigorous mathematical results have justified this approach and have identified some common features that arise for certain classes of integrable models.  In this work we present a thorough numerical investigation that shows that several of these core features extend to a paradigmatic linear-in-parameters model, the Lotka-Volterra (LV) system, which we consider in the conservative case as well as under the addition of terms that perturb the system away from this regime.
 A central construct for this analysis is a concise representation of parameter features in the data space that we call the $P_n$-diagram, which is particularly useful for visualization of results for low-dimensional (small $n$) systems.
Our work also exposes some new  properties related to non-uniqueness that arise for these LV systems, with non-uniqueness manifesting as a multi-layered structure in the associated $P_2$-diagrams.
\end{abstract}

\submitto{\IP}
\maketitle

\section{Introduction}
A fundamental step in computational modeling is to fit a mathematical model of a physical system to data collected from observations of one or more instantiations of that system.  
The mathematical model may be derived from physical principles or may be a phenomenological representation designed to capture system dynamics.  Either way, if the model evolves continuously in time and the rate of change of the state variables is a function of the state, the model can be written as a system of ordinary differential equations
\begin{equation}
    \label{eq:genmodel}
    d\mathbf{x}(t)/dt = \mf(\mathbf{x}(t),\mathbf{p})
\end{equation}
with variables or states $\mathbf{x} = (x_1, \ldots, x_n)$ and parameter vector $\mathbf{p}$ in a permissible set $\Omega$. 
In this setting, the fitting process takes the form of an \emph{inverse problem} in which one aims to solve for the components of $\mathbf{p}$ from measurements of either the trajectory $ \mathbf{x}(t)$ or some known function of the trajectory (i.e., the output), $ \mathbf{y}(t)=\mathbf{h}(\mathbf{x}(t))$.


We here consider the inverse problem for
situations in which we have discrete data points derived from a single experiment. 
%
If equation (\ref{eq:genmodel}) were a perfect,  error-free representation of the modeled physical system and all variables were observable, then this situation would correspond to using data points $P_j$ with coordinates $\{\mathbf{x}_j= \mathbf{x}(t_j) = \mathbf{\phi}(t_j;\mathbf{p^*}) \}$ for some times $\{ t_j \}$, where $\mathbf{\phi}(t,\mathbf{p^*})$ is a trajectory generated by (\ref{eq:genmodel}) for some choice $\mathbf{p^*}$ of parameter values.
In fact, the restriction to data from a single experiment arises in modeling many real-world systems for which experiments are non-repeatable, including living systems in the lab or the clinic with heterogeneity across individuals as well as climate and other systems that cannot be set to specific initial conditions by an experimentalist.

Real data are often subject to experimental noise, which is why much of the literature on the subject of parameter estimation employs statistical approaches, where the inverse problem is formulated probabilistically and the likelihood of the estimated parameters is of primary interest\cite{tarantola2005inverse, calvetti2007introduction,stuart2010inverse, smith2013uncertainty, calvetti2018inverse}. As a complement to those approaches, we here focus on more fundamental aspects of the inverse problem that show up even in the absence of measurement error. Among these are problems of existence, uniqueness (or identifiability), and robustness (or sensitivity of parameters on data) of solutions.
That is, if we fix $\mathbf{f}$ and obtain a set of measured data points, then it is possible that no $\mathbf{p} \in \Omega$ exists for which there is a trajectory $\mathbf{\phi}(t;\mathbf{p})$ of (\ref{eq:genmodel}) passing through the data points.   
In addition, there may be more than one choice of $\mathbf{p}$ that gives rise to such a $\mathbf{\phi}(t;\mathbf{p})$.  This situation can occur due to non-identifiability, which refers to a situation in which the measured data depends only on combinations of components of $\mathbf{p}$ that remain constant over sets of $\mathbf{p}$ values in $\Omega$, or due to properties of $\mathbf{f}$.
Finally, although the trajectories (and hence data) depend differentiably on $\mathbf{p}$, the inverse need not be true, i.e., the data may not depend differentiably or even continuously on the data, which makes the inverse problem ill-conditioned. As a result, small changes in data 
will have great influence on inferred parameters and our understanding of the system. As discussed in \cite{swigon2019importance}, all of these complications need to be understood before we even begin to try to incorporate random effects and measurement noise.

Although these challenges are daunting, mathematical theory gives us tools to address them systematically by exploring the details of the relationship between parameters and data. In particular, our approach is to combine classical results on existence and uniqueness of solutions of ODEs and their differentiable dependence on parameters \cite{hartman2002ordinary,teschl2012ordinary} with basic results from singularity theory on invertibility of nonlinear maps \cite{arnold1981singularity, krantz2002implicit} to provide necessary and/or sufficient \emph{conditions on the data} that provide desired qualitative charactristics of the inverse problem. The emphasis on conditions on the data distinguishes our work from related work on identifiability, observability, and controlability, which has focused on conditions satisfied by the parameters of the system (see, e.g., \cite{kalman1963mathematical,kou1973observability,aeyels1981generic,griffith1971observability, Miao2011}).

Because of the generality of equation (\ref{eq:genmodel}) and the associated inverse problem, we do not expect a one-size-fits-all theory that encompasses all such scenarios.
To make headway in the case of inverse problem from a single trajectory, we began to develop mathematical theory in the setting of linear systems with specific assumptions on the data available \cite{stanhope2014,stanhope2017}.  From that starting point, we achieved new analytical results that are summarized in Section \ref{linearresults} below with new results on invariance of $P_n$ diagrams described by Theorem \ref{linearP2invariance}. Those were later generalized to affine systems \cite{duan2020identification} (see Section \ref{affineresults} with new theorem \ref{affineP3invariance}). The conjecture underlying this systematic approach is that these results, or at least some of their features, can be generalized in directions that include allowing specific forms of nonlinearity in $\mathbf{f}$, and in this work, we desribe one such generalization by considering how our results extend to an important class of nonlinear systems, the linear-in-parameters (LIP) systems, and what new features may emerge in this nonlinear setting.  This study provides a natural opportunity to collect results from across several model classes and integrate them in a way that, we hope, will be useful to guide continued work in related directions.

\section{Lessons learned from integrable systems}

Existence of explicit solutions of the initial value problem for \eref{eq:genmodel} reduces investigations of the inverse problem to the study of inverses of nonlinear algebraic functions \cite{krantz2002implicit}. A thorough analysis of linear, affine, and matrix Riccati systems reveals a picture of a complex, but understandable, dependence of the existence, uniqueness, and properties of the inverse problem solution on the data, which is outlined in subsections below. In each case we find it revealing to display the results of this dependence schematically using what we henceforth call the \emph{$P_n$-diagram} (with $n$ replaced by an appropriate integer). This diagram depicts the dependence of the existence, uniqueness, and other qualitative features of the solution of the inverse problem as a function of the last data point $P_n$ (with coordinates $\mathbf{x}_n$) when the first $n$ data points $P_0, P_1, ..., P_{n-1}$ (with coordinates $\mathbf{x}_0,...,\mathbf{x}_{n-1}$) are held fixed.

\subsection{Linear systems}\label{linearresults}
Linear systems provide a suitable starting point for theoretical investigations of parameter identifiability and estimation because their solutions can be readily expressed in closed form.  
These systems take the form 
\begin{equation} \label{eq:linear}
         d\mathbf{x(t)}/dt = A\mathbf{x}(t),
\end{equation}
with dependent variable $\mathbf{x} \in \mathbb{R}^n$, parameters $A \in \mathbb{R}^{n \times n}$, and initial condition $\mathbf{x}(0) = \mathbf{b} \in \mathbb{R}^n$ for a natural number $n$. 
Despite the linearity of the ODE in (\ref{eq:linear}), the relation between the parameters and the data 
is nonlinear; thus, we can already start to derive non-trivial insights into the inverse problem in this dynamically simple setting.

For a unique solution to the inverse problem to exist, it is intuitively natural that we need to have at least as many pieces of information available as we have parameter values in the system.  For system (\ref{eq:linear}), there are $(n+1)n$ such parameters consisting of the entries of $A$ and $\mathbf{b}$.  For simplicity we assume that all components of $\mathbf{x}$ can be observed; and with that assumption our information requirement amounts to the need for a set of at least $n+1$ data points, $d := \{ P_j  : j = 0, 1, \ldots, n \}$ with coordinates $\mathbf{x}_j \in \Rn$. Let us focus for now on the case when we have exactly $n+1$ data points, collected at equally spaced times $t_j$ such that $t_{j+1}-t_j = \Delta t$ is constant; moreover, for simplicity, we will take $\Delta t = 1$ and translate time as needed to take $\mathbf{b}=\mathbf{x}_0$. (Extensions to non-equally spaced points and a subset of variables being observed can be found in \cite{stanhope2017}.)  
 
With this set-up, we can define forward and inverse problems and associated mappings.  Use ${\cal D}$ to denote the data space, consisting of sets of collections of $(n+1)$ points (i.e., $(n+1)$-tuples) in $\mathbb{R}^n$, and ${\cal P}$ to denote the parameter space, comprising all sets of two elements such that one is a matrix in $\mathbb{R}^{n \times n}$ and the other is a vector in $\mathbb{R}^n$.
The forward problem is the problem of finding the mapping $F: {\cal P} \to {\cal D}$; this problem can also be referred to as solving the initial value problem (\ref{eq:linear}), with $F$ as the solution map.  The inverse problem is the problem of inverting $F$, such that from a given data set $d \in {\cal D}$, we can determine the $(A,\mathbf{b}) \in {\cal P}$ for which $F(A,\mathbf{b})=d$; that is, we can solve $F^{-1}(d)=(A,\mathbf{b})$.  Solving the inverse problem is sometimes referred to as parameter identification or parameter estimation; the latter, however, is often used to refer to a more general scenario in which the information obtained from $d$ may be estimates, rather than exact values, of $A$ and $\mathbf{b}$.
 
In this setting, we have obtained four \emph{fundamental results} on the inverse problem for system (\ref{eq:linear}), which we concisely review here.
\emph{First}, a necessary but not sufficient condition for the uniqueness of a solution $A$, with $\mathbf{b}=x_0$, is that $d$ is a linearly independent set \cite{stanhope2014}.
\emph{Second}, if $d$ is linearly independent, then necessary and sufficient conditions for existence of an inverse problem solution and for the uniqueness of the solution when it exists can be expressed in terms of the eigenvalues and Jordan block structure of the matrix $ \Phi = X_1 (X_0)^{-1}$,
where $ X_k, \; k \in \{ 0, 1 \},$ is the square matrix with columns  $[ \mathbf{x}_k | \mathbf{x}_{k+1} | \ldots | \mathbf{x}_{k+n-1} ]$,
and the invertibility of $X_0$ follows from the linear independence of $d$ \cite{stanhope2017}.
\emph{Third}, under slightly stronger conditions on the spectrum of $\Phi$, small perturbations to all elements of $d$ do not change the qualitative nature of the inverse problem solution; rather, existence or non-existence, uniqueness, and even various other properties of the solution matrix are preserved \cite{stanhope2017}:

\begin{theorem}
 \label{thm:opensets}
 Suppose that for $d \in {\cal D}$, $X_0$ is invertible and $\Phi = X_1 X_0^{-1}$ has
 \begin{description}
 \item{(a)} only positive real or complex eigenvalues, 
 \item{(b)} distinct positive real eigenvalues, or
 \item{(c)} at least one negative eigenvalue of odd multiplicity.
 \end{description}
 Then there exists an open set $U \subset {\cal D}$ containing $d$ such that for every $d^* \in U$, the mapping $F^{-1}(d^*)$ has, respectively,
 \begin{description}
 \item{(a)} a real solution,
 \item{(b)} a unique real solution, or
 \item{(c)} no real solution.
 \end{description}
 \end{theorem}
 
\emph{Fourth}, analytical lower and upper bounds relating to the size of $U$ have been computed \cite{stanhope2017}.  Specifically, for a given $d$ and a property of interest, which may be existence, non-existence, or uniqueness of an inverse problem solution, or even stability of the system, perturbations to $d$ that are smaller than the lower bound are guaranteed to maintain the solution property for which the bound is derived, while for any size larger than the upper bound, there is certain to exist a perturbation of that size for which the inverse solution does not have that property. 
Since each of these properties depends on conditions on the eigenvalues of $\Phi$ in the complex plane, it is not surprising that the lower bound grows with the distance of these eigenvalues from the boundary of the region where the conditions hold, while it shrinks in proportion to how strongly perturbations to $\Phi$ impact its spectrum. The latter property might naturally be thought on as relating to the condition number of $\Phi$, but it does so in a way that depends on the location of the data points.  The upper bound, in turn, relates to the perturbation needed to guarantee a breakdown of the eigenvalue properties required for the property of interest. 

\begin{figure}
\centering
\includegraphics[width=0.8\textwidth]{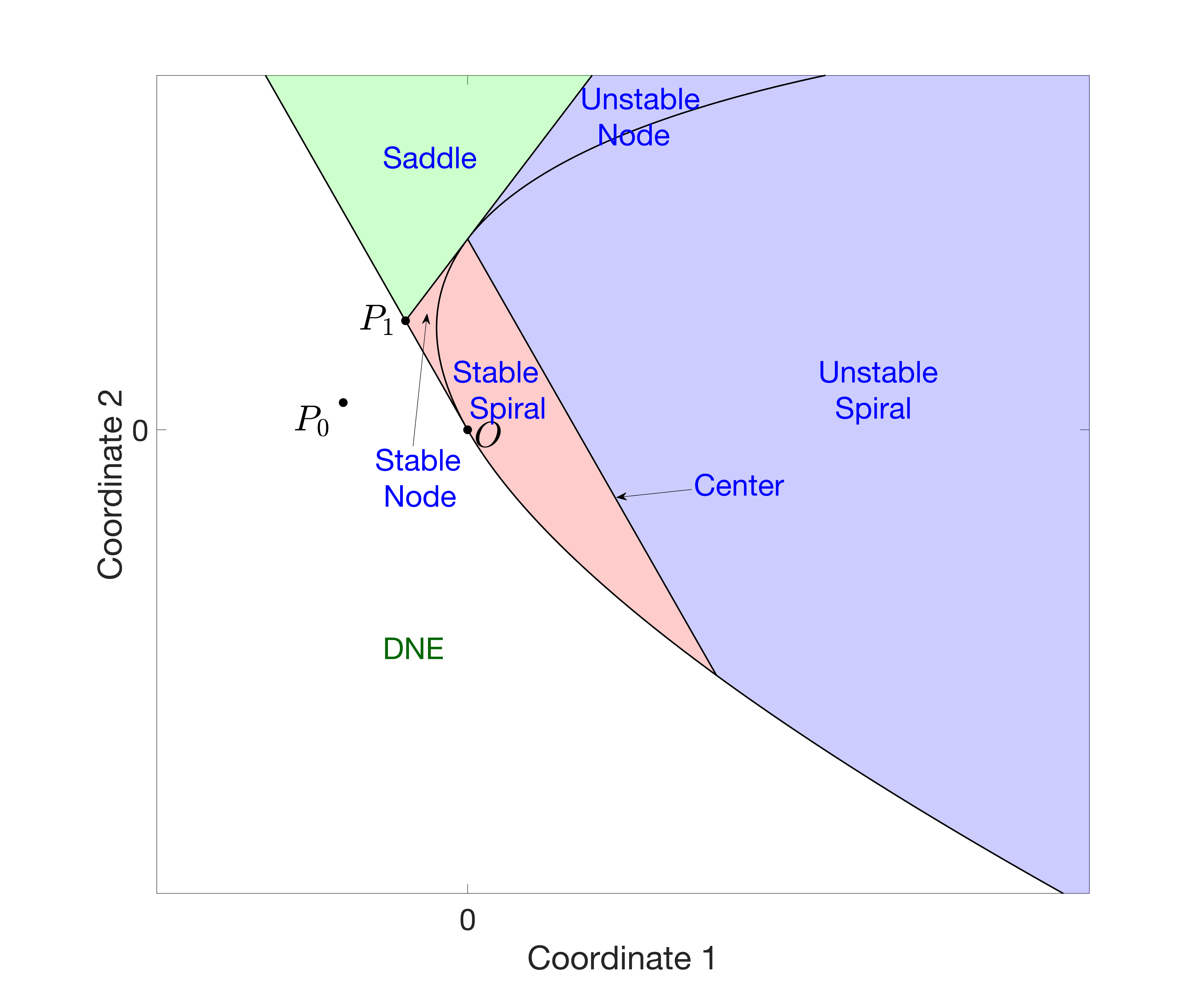}
\caption{The $P_2$-diagram for linear system \eref{eq:linear} with $n=2$.  This image depicts the regions of location of $P_2$ for which the inverse problem for system \eref{eq:linear} either does not have a solution (DNE) or has a solution such that the origin O has the specified stability properties when the locations of data points $P_0,P_1$ are fixed as indicated.}  
\label{figure:P2diagramlinear}           
\end{figure}

These results generalize to other properties of the inverse problem solution specific to linear systems, as long as those properties can be reduced to conditions on the eigenvalues of the parameter matrix $A$.  
For example, if a given $d$ has a unique inverse problem solution $F^{-1}(d) = (A,\mathbf{b})$ such that the origin is a stable node for the linear system $\mathbf{x}' = A\mathbf{x}$, then as long as some non-degeneracy conditions on the eigenvalues of $\Phi$ hold, the origin will remain a stable node for $F^{-1}(d^*)$ for all $d^*$ in a neighborhood of $d$, the size of which we can estimate.
More importantly, Theorem \ref{thm:opensets} and its extensions to other properties of the dynamical system reveal nice features of the regions in data space that correspond to the solution of the inverse problem. 
Specifically, once data points $\hat{d} := \{ P_0, \ldots, P_{n-1} \}$ are fixed, the $n$-dimensional space parameterized by the coordinates $\mathbf{x}_n$ consists of a collection of open regions such that if two possible choices of $P_{n_1}, P_{n_2}$ are made from the same region, then the inverse problem solution will have the same properties for $d_1 := \{ \hat{d}, P_{n_1} \}$ and for $d_2:= \{ \hat{d}, P_{n_2} \}$.  This is exemplified in 2 dimensions by the $P_2$-diagram depicted in Figure \ref{figure:P2diagramlinear}. This diagram makes clear that the boundaries of these regions are smooth curves given by simple algebraic expressions\cite{stanhope2017}. Furthermore, this diagram is universal, in that although it depends on the positions of the data points $P_0,P_1$, that dependence is defined by simple affine transformation of the 2D coordinate system given by stretching and rotation about the origin, with the possible addition of a mirror image transformation. Formally, this is a result of the invariance of solutions of the dynamical system \eref{eq:linear} to transformations as follows:

\begin{theorem}\label{linearP2invariance}
Suppose that a dynamical system \eref{eq:linear} has, for data points $P_j, j=0,1,..,n$ with coordinates $\mathbf{x}_j$, an inverse problem solution given by the parameter matrix $A$. Then for any invertible matrix $S$ the same system has an inverse problem solution for data point coordinates $\bar{\mathbf{x}}_j = S\mathbf{x}_j $ with $j=0,1,..,n$. This solution is given by the parameter matrix $\bar{A}=SAS^{-1}$, which has the same eigenstructure as $A$.
\end{theorem}

\begin{proof}
Clearly, if there is a solution $A$ of the inverse problem with given data point coordinates $\mathbf{x}_j, j=0,1,..,n$, then the system \eref{eq:linear} with parameter matrix $A$ has a trajectory $\mathbf{x}(t)$ that passes through the data, i.e., $\mathbf{x}_j = \mathbf{x}(j) $ for $j=0,1,..,n$.
Let $S = \bar{X}_0X_0^{-1}$ where $X_0$ is as above and $\bar{X}_0 = [ \bar{\mathbf{x}}_0 | \bar{\mathbf{x}}_{1} | \ldots | \bar{\mathbf{x}}_{n-1} ]$.  Clearly, $\bar{\mathbf{x}}(t) = S\mathbf{x}(t)$ solves the system \eref{eq:linear} with matrix $SAS^{-1}$ and such a solution passes through the data $\bar{\mathbf{x}}_j = \bar{\mathbf{x}}(j) $ for $j=0,1,..,n$.
\end{proof}

Theorem \ref{linearP2invariance} makes clear that the $P_n$-diagram is determined by the position of the first $n$ data points labeled $P_0, P_1,...,P_{n-1}$. If those points are changed, there is a one-to-one transformation between the points of the original and the transformed diagram, realized as multiplication by the matrix $S$. The preservation of eigenstructure of the parameter matrix upon transformation results in preservation of the stability properties of the system and hence the nature of the domains in the $P_n$-diagram and their boundaries.

 
 \subsection{Affine systems}\label{affineresults}
 A first, relatively gentle extension of linear systems is provided by consideration of affine systems, which take the form
 \begin{equation}
     \label{eq:affine}
     d\mathbf{x(t)}/dt = A\mathbf{x}(t) + \mathbf{c}
 \end{equation}
 with states $\mathbf{x} \in \mathbb{R}^n$ and parameters $A \in \mathbb{R}^{n \times n}, \mathbf{c} \in \mathbb{R}^n$, and initial condition $\mathbf{x}(0) = \mathbf{b} \in \mathbb{R}^n$.
 Affine systems, like linear systems, can be solved explicitly.  In fact, if $A$ is invertible with $\mathbf{c} \in \mbox{range}(A)$, then system (\ref{eq:affine}) has a unique equilibrium point at $\mathbf{x}^* := -A^{-1} \mathbf{c}$ and can be transformed into a linear system in $\my$ for $\my = \mx-\mx^{*}$.
 
 For general affine systems (\ref{eq:affine}), even without the assumption of invertibility of $A$, many of the identifiability properties shown to hold for linear systems (\ref{eq:linear}) still carry over \cite{duan2020identification}.
 The key step in establishing this generalization is to define the augmented linear system
 \begin{equation}
     \label{eq:aug}
 d\mz(t)/dt =  B\mz := \left[ \begin{array}{cc} A &  \mc \\ 0 & 0 \end{array} \right] \mz, \; \; \mz(0) = \left[ \begin{array}{c} \mb \\ 1 \end{array}
 \right]
 \end{equation}
 for $ \mz(t) := [ \, \mx(t) \; 1 \, ]^T  $, with solution matrix $\Psi$.  If we define $Z_0, Z_1$ analogously to $X_0, X_1$, using $\mz_i = [ \, \mx_i \; 1 \, ]^T$ in place of $\mx_i$ for each $i$, and if $Z_0$ is invertible, then we can also write
 $\Psi = Z_1 Z_0^{-1}$, which is useful for studying the inverse problem for (\ref{eq:affine}). 
 
\begin{figure}
\centering
\includegraphics[width=0.8\textwidth]{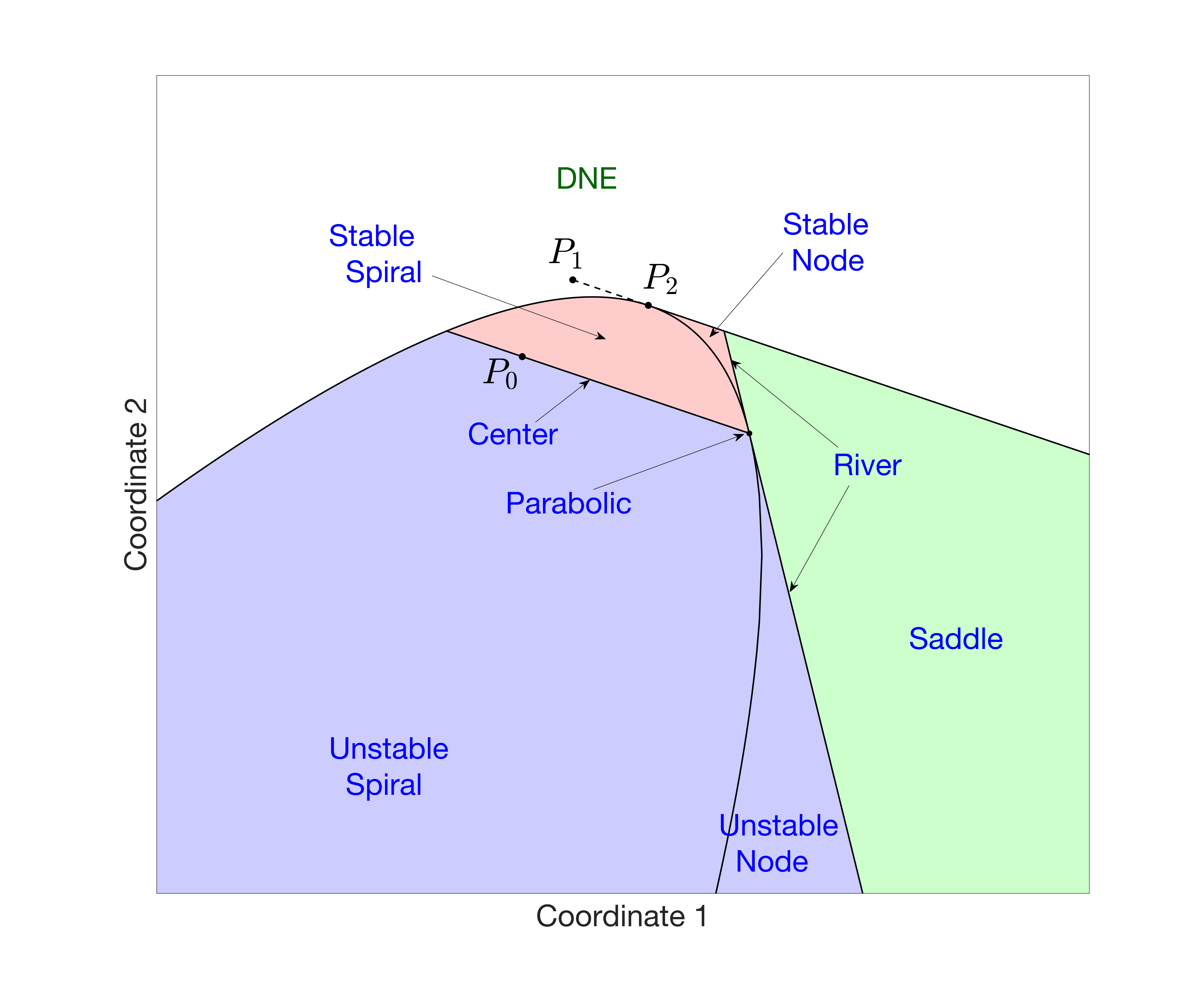}
\caption{The $P_3$-diagram for affine system \eref{eq:affine} with $n=2$, which depicts the regions of location of $P_3$ for which the inverse problem has a solution with specified stability properties when the location of data points $P_0,P_1,P_2$ is fixed as indicated. Note that in the river and parabolic regimes, the system has no fixed point, which cannot occur in the linear case.}  
\label{figure:P3diagramaffine}           
\end{figure}

As in the linear case, we  find that linear independence of data points is helpful for the uniqueness of a solution to the inverse problem for (\ref{eq:affine}) and that the properties of the solution relate to the spectrum of a matrix that can be derived from the data.  In fact, in the affine case, the matrix $\Psi$ relates to the solution matrix $\Phi$ for (\ref{eq:linear}) and hence, not surprisingly, properties of the spectrum of $\Phi$ also impact the solution of the inverse problem for (\ref{eq:affine}).  Moreover, the existence of open regions in data space where inverse problem features persist also extends to the affine case, with some modifications to the analytical lower bounds on the sizes of regions where persistence occurs \cite{duan2020identification}. In $n$-dimensional state space such regions make up the $P_{n+1}$-diagram for the inverse problem, which for $n=2$ is depicted in Figure 3(a) of \cite{duan2020identification} and reproduced here in Figure \ref{figure:P3diagramaffine}. Similarly to the linear case, this diagram depends on data points $P_0,P_1,P_2$ only via a simple transformation of the state space consisting of stretching, rotation, and mirror inversion, combined with translation. Formally this is represented by the following result:

\begin{theorem}\label{affineP3invariance}
Suppose that a dynamical system \eref{eq:affine} has, for data points $P_j, j=0,1,..,n+1$  with coordinates $\mathbf{x}_j$, a solution of the inverse problem given by $(A,\mathbf{c})$. Then for any invertible matrix $S$ and vector $\mathbf{r}$ the system \eref{eq:affine} has an inverse problem solution for data point coordinates $\bar{\mathbf{x}}_j = S\mathbf{x}_j + \mathbf{r} $ with  $j=0,1,..,n+1$. This solution is given by the ordered pair  $(\bar{A},\bar{c}) = (SAS^{-1},S(\mathbf{c}-AS^{-1}\mathbf{c}))$.
\end{theorem}

\begin{proof}
Clearly, if there is a solution $(A,\mathbf{c})$ of the inverse problem with given data $\mathbf{x}_j, j=0,1,..,n$, then the system \eref{eq:affine} with parameters $(A,\mathbf{c})$ has a trajectory $\mathbf{x}(t)$ that passes through the data, i.e., $\mathbf{x}_j = \mathbf{x}(j) $ for $j=0,1,..,n+1$.
Let $R = \bar{Z}_0Z_0^{-1}$ where $Z_0$ is as above and $\bar{Z}_0 = [ \bar{\mathbf{z}}_0 | \bar{\mathbf{z}}_{1} | \ldots | \bar{\mathbf{z}}_{n} ]$ with $\bar{\mz}_j = [ \, \bar{\mx}_j \; 1 \, ]^T$ and let matrix $S$ and vector $\mathbf{r}$ be such that
\[
R = \left[ \begin{array}{cc} S &  \mathbf{r} \\ 0 & 1 \end{array} \right].
\]
It can be easily verified that $\bar{\mathbf{x}}(t) = S\mathbf{x}(t) + \mathbf{r}$ solves the system \eref{eq:affine} with parameters $(SAS^{-1},S(\mathbf{c}-AS^{-1}\mathbf{r}))$ and that such a solution passes through the data as $\bar{\mathbf{x}}_j = \bar{\mathbf{x}}(j) $ for $j=0,1,..,n$.
\end{proof}

Theorem \ref{affineP3invariance} makes clear that in $n$-dimensional state space the $P_{n+1}$-diagram for affine systems is determined by the position of the first $n+1$ data points  $P_0,P_1,...,P_n$. If those points are changed, there is a one-to-one transformation between the points of the original and the transformed diagram, realized with the matrix $S$ and translation $\mathbf{r}$. The preservation of eigenstructure of the parameter matrix upon transformation results in preservation of the stability properties of the system and hence the nature of the domains in the $P_{n+1}$-diagram and their boundaries.

\subsection{Nonsymmetric matrix Riccati systems}
An interesting extension of the above results arises from considering a nonlinear dynamical system provided by the nonsymmetric matrix Riccati differential equation (RDE):
  \begin{equation}
     \label{eq:RDE}
 dX(t)/dt =  -X(t)BX(t) - X(t)A + DX(t) + C
 \end{equation}
with states $X \in \Rmn$, parameter matrices $A \in \Rnn$, $B \in \Rnm$, $C \in \Rmn$, $D \in \Rmm$, and initial condition $X(0) = X_0 \in \Rmn$.
Classical theory of Radon tells us that, like the affine system, any solution of this system is locally equivalent to a solution of a linear ODE system with block coefficient matrix $M = [A\ B; C\ D]$ (see \cite{freiling2002survey}). This result can be used to solve the inverse problem for the system \eref{eq:RDE}, i.e., to determine the parameter matrices $A,B,C,D$ from data $d=\{X_0,X_1,...,X_p\}$ for a certain $p$ \cite{swigonparameter}. Given that there are $(m+n)^2$ parameters and $m+n$ initial conditions in the problem, one expects that unique determination of the parameters can be achieved with $p=\lceil(m+n+1)(m+n)/(mn)\rceil$ data points, with the exception of a structural unidentifiability in the system where changing $A$ into $A+\alpha I$ and $D$ into $D+\alpha I$ does not change the dynamics.
Although a detailed analysis of this case is yet to be done, the availability of explicit solutions and formulas for the parameter matrices in terms of the data makes it straightforward to formulate conditions for existence and uniqueness of a solution to the inverse problem for \eref{eq:RDE} and hence to generalize the fundamental results first observed for linear systems. As an example, a straightforward extension of Theorems \ref{linearP2invariance} and \ref{affineP3invariance} can be formulated to show that the $P_p$-diagram is invariant with respect to transformations of the data point coordinates of the form $\bar{X}_j = S(X_j+R)Q$ where $S\in \Rmm$ and $Q\in\Rnn$ are invertible matrices and $R\in \Rmn$.

 \section{Linear-in-parameters models}
 The extension of fundamental properties of the inverse problem from linear to affine systems raises the intriguing possibility that they may hold for more general classes of nonlinear systems.
 The class that we will consider in the remainder of this work are the linear-in-parameters (LIP) systems, which take the form
 \begin{equation}
     \label{eq:LIP}
     d\mx(t)/dt = B\mf(\mx(t))
 \end{equation}
 for $\mx \in \mathbb{R}^n$, $B \in \mathbb{R}^{n \times m}$, and $\mf:\mathbb{R}^n \to \mathbb{R}^m$.
 
Previous numerical work \cite{duan2020identification} showed that affine approximation of this system can be highly effective for estimating the inverse problem solution $B$ for (\ref{eq:LIP}) from a discrete data set.
Indeed, for continuously differentiable $\mf(\mx)$, linearization of (\ref{eq:LIP})  around any point $\mx^*$ in the domain of $\mf$ yields
\[
\mf(\mx) = \mf(\mx^*) + \nabla \mf(\mx^*)(\mx-\mx^*) + R(\mx,\mx^*)
\]
for $||R(\mx,\mx^*)|| \sim \Or(||\mx-\mx^*||^2)$.
Thus, system (\ref{eq:LIP})  is well approximated in a neighborhood of $\mx^*$ by the affine system $d\mathbf{w}/dt = A\mathbf{w}(t)+\mathbf{c}$ for $A = B \nabla \mf(\mx^*)$ and $\mathbf{c} = B(\mf(\mx^*)-\nabla \mf(\mx^*)\mx^*)$.   This result, however, does not provide any insight into the general properties of inverse problem solutions for such LIP systems.

To date only the first fundamental result has been generalized as follows \cite{stanhope2014}: for a given trajectory $\gamma(B) := \{\mx(t;B) | t \in \mathbb{R} \}$ of system (\ref{eq:LIP}), the uniqueness of $B$ as a solution to the inverse problem holds if and only if $\{ \mf(\mx(t;B) | t \in \mathbb{R} \}$ is not confined to a proper subspace of $\mathbb{R}^m$. In terms of discrete data, the necessary (but not sufficient) condition for uniqueness of $B$ is therefore that the images $\{\mf(\mx_0),\mf(\mx_1),...,  \mf(\mx_p)\}$ are linearly independent.

The lack of explicit solutions of the forward problem, and thereby a lack of explicit or implicit functions that describe the dependence of such solutions on the parameters of the system, hinders any systematic description of regions in data space that correspond to existence, uniqueness, or various qualitative aspects of the solutions of the inverse problem. We therefore turn to numerical studies of example nonlinear systems to gauge the range of validity of the fundamental results observed for linear and affine systems and to identify departures from those results. We use boundary value problem solvers combined with continuation methods to  numerically compute the boundaries between existence/nonexistence regions and various types of dynamics in data space and describe the results below. Because we have seen its utility for summarizing much information about the solvability and solutions of the inverse problem in other cases, we are particularly interested in finding the $P_n$-diagram for LIP systems and observing both how it transforms with changes in data point coordinates and how it differs qualitatively from the diagrams found for linear and affine systems.

\section{A Lotka-Volterra system}
\label{section:mainresult}
As we have discussed, a key idea from the study, in linear and affine systems, of the inverse problem of determining $(n+1)n$ unknown parameters from $n+1$ discrete, linearly independent data points in $\mathbb{R}^n$ is that if we treat $n$ of these data points as fixed, then we can partition $\mathbb{R}^n$ into simply bounded regions such that the region in which the $(n+1)$st data point lies determines the properties of the inverse problem solution.
To test the extent to which this idea carries over to LIP systems, we performed a thorough numerical study of this approach to the inverse problem for discrete data in a specific such system, namely the famous Lotka-Volterra (LV) equations, which have a rich history, wide range of applications, and a well developed theory \cite{wangersky1978lotka,takeuchi1996global,anisiu2014lotka}.

We consider the LV equations in the form
\begin{equation}\label{LV3}             
\left\{
             \begin{array}{lr}
             \dot{x}=\aq x+\bq xy, &  \\
             \dot{y}=\bw xy+\aw y, &  \\
             x(0)=x_0, &\\
             y(0)=y_0. &
             \end{array}
\right.
\end{equation}
Advantages of this choice are the low-dimensionality of the model, with $n=2$, and the fact that the model has potential biological relevance over a wide range of parameter sign structures. The model is not integrable except for special parameter combinations \cite{brenig1988complete} but one can construct power-series approximations of solutions \cite{olek1994accurate}. We focus on positive solutions of (\ref{LV3}), since $x$ and $y$ generally denote either the numbers or the densities of species in two interacting populations, the first quadrant is invariant under the flow of (\ref{LV3}), and the results for any other quadrant (likewise invariant) can be related to those in the first quadrant by a proper adjustment of the signs of parameters.  Note that system (\ref{LV3}) always has a critical point at the origin but need not have a critical point in the interior of the positive quadrant.
 
One can write the LV system in the LIP form (\ref{eq:LIP}) as follows:
\begin{equation}\label{LIP}        
\left\{
             \begin{array}{lr}
             \dot{\mx}=A\mf(\mx), &  \\
             \mx(0)=\mb,&
             \end{array}
\right.
\end{equation}
where
\begin{equation}
\label{LV3alt}
\mx=\left[
\begin{array}{c}
x\\
y
\end{array}
\right],
A=\left[
\begin{array}{ccc}
\aq & \bq & 0\\
0 & \bw & \aw
\end{array}
\right],
\mf(\mx)=\left[
\begin{array}{c}
x\\
xy\\
y
\end{array}
\right],\rm{ and } \; 
\mb=\left[
\begin{array}{c}
x_0\\
y_0
\end{array}
\right].
\end{equation} 
For brevity of notation we will use the matrix $A$ to represent the parameters $(\aq,\bq,\bw,\aw)$ of the system \eref{LV3}, implicitly assuming that in any such matrix $A$ the entries $a_{13}$ and $a_{21}$ are equal to zero and never take any other values. We denote by $\sigma_A$ the \emph{signature of the system}, i.e., 
\begin{equation}
\sigma_A=[\rm{sgn}(\aq)\quad  \rm{sgn}(\bq)\quad \rm{sgn}(\bw)\quad \rm{sgn}(\aw)],
\end{equation}
which determines the type of interactions described by the model.
In the classical theory, positive values of both $\bq$ and $\bw$ correspond to a cooperative interaction, negative values of both correspond to a competitive interaction, while the system is classified as a predator-prey system when $\bq$ and $\bw$ have distinct signs.  

We will primarily consider a data set $d = \{P_0, P_1, P_2\}$  where $P_j, j=0,1,2,$ has coordinates $(x_j,y_j)$ such that $0 < x_0 < x_1$ and $0 < y_0 < y_1$. As before, will assume that $P_0$ and $P_1$ are fixed and that $P_2$ can lie anywhere in the first quadrant of the $(x,y)$ plane.
Within this domain we ask the following questions:
\begin{enumerate}
\item Existence: What is the set of values of $P_2$ for which there exists some $A$ (i.e., parameters $(\aq,\bq,\bw,\aw)$), such that the system defined by \eref{LV3} or equivalently \eref{LIP}, \eref{LV3alt} has a trajectory $\phi(t;A)$ with $P_j=\phi(j,A)$, $j\in\{0,1,2\}$?

\item Uniqueness: What is the set of values of $P_2$ for which the parameter matrix $A$ that solves the inverse problem is unique?

\item Parameter properties:  What is the set of values of $P_2$ for which the parameter matrix $A$ has specific signs of its entries, $\sigma_A$ (corresponding to specific types of behavior)? 
\end{enumerate}

\begin{figure}
\centering  
\includegraphics[width=1\textwidth]{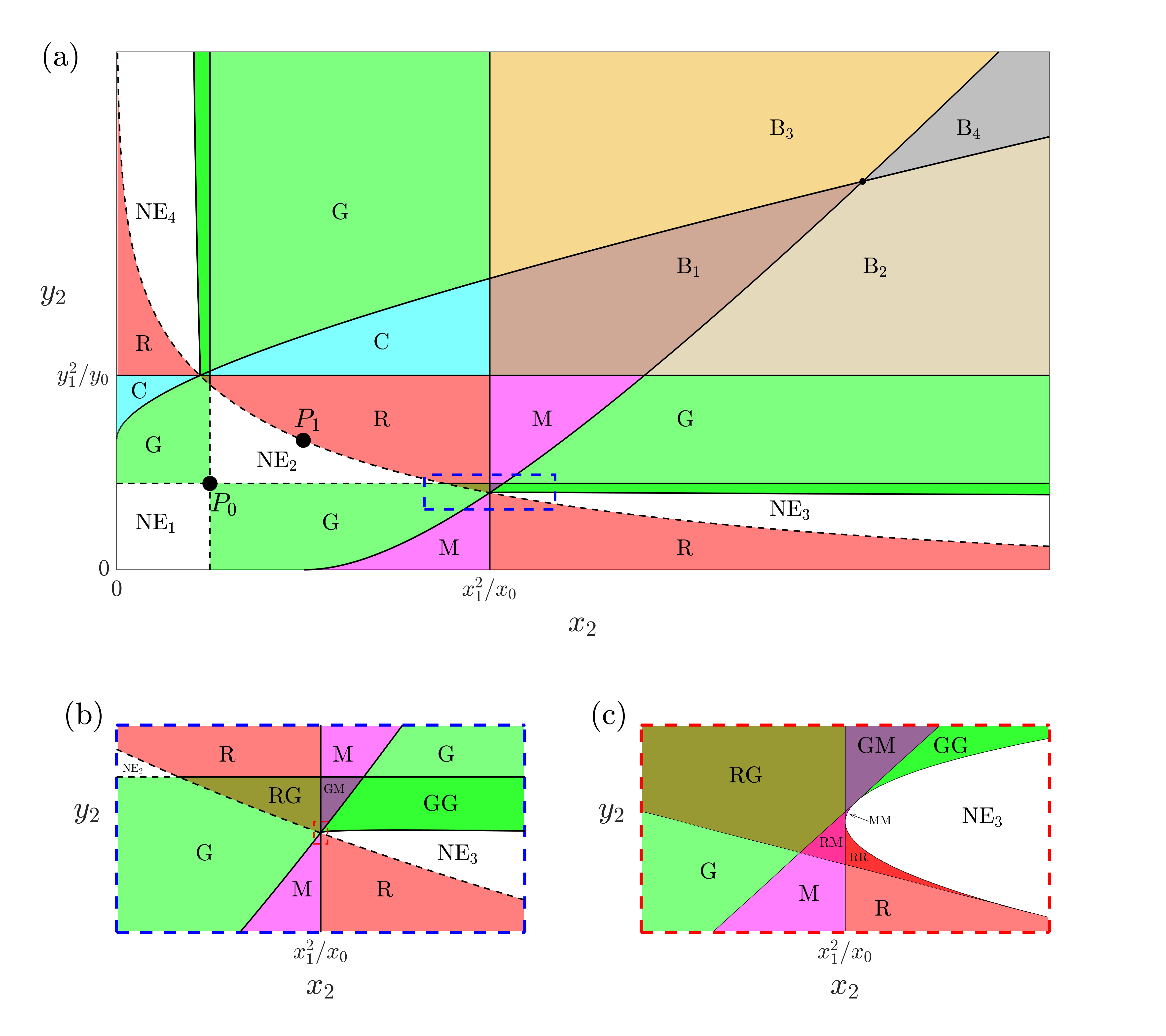}
\caption{The $P_2$-diagram for Lotka-Volterra system \eref{LV3}, which depicts the regions $\mathcal{R}_\Omega$ of $P_2$ points for which the inverse problem has a specified signature $\sigma_A$ (only the subscripts $\Omega$ of the  regions are shown).  Panels (b) and (c) are enlargements of the interior dashed rectangles in panels (a) (blue dashed) and (b) (red dashed), respectively. The data points $P_0$ and $P_1$ are fixed at the labeled positions. Solid curves represent region boundaries for which inverse problem solutions exist, while dashed curves boundaries where solutions do not exist. See text and Table \ref{table:P2diagram} for more detailed descriptions of the regions.}
\label{figure:P2diagram}           
\end{figure}

Via extensive use of numerics, described in Section 1 of Supplementary Materials, we find that the desirable properties obtained in linear and affine systems carry over to this LIP system; specifically, we obtain clearly defined, simply bounded regions such that the properties of the inverse problem solution (or lack thereof) are the same for all choices of data points $P_2$ within the same region. The numerical results can be interpreted in two ways: (a) as a manifold $\mathcal{M}$ in 6-dimensional space with coordinates $(x_2,y_2,\aq,\bq,\bw,\aw)$ comprised of points for which parameter matrix $A$ constructed from parameters $\aq,\bq,\bw,\aw$  solves the inverse problem with $P_2$ located at $(x_2,y_2)$, or (b) as a projection of the manifold $\mathcal{M}$ onto the 2-dimensional space with coordinates $(x_2,y_2)$, shown in Figure \ref{figure:P2diagram} and described in Table \ref{table:P2diagram}, which, in analogy to the earlier sections, we call the \emph{$P_2$-diagram}. The $P_2$-diagram contains a large amount of information and depicts schematically the answers to questions 1-3 above; that is, it displays the  locations of points $P_2$ for which we have evidence about the existence, uniqueness and properties of the inverse problem solution.

\begin{table}[ht]
\centering
\caption{Properties of systems with $P_2 \in \mathcal{R}_\Omega$ of Figure \ref{figure:P2diagram} }
\label{table:P2diagram}           
\begin{tabular}{|p{1.6cm}|p{2.3cm}|p{5.4cm}|p{3.4cm}|} 
\hline 
\ct $\Omega$ & \ct Number of solutions 
&  $\sigma_A$ & Dynamics type\\
\hline 
\ct $\rm{R}$ & \ct $\geq$1 &$[+--+]$ & competitive\\
\hline 
\ct $\rm{G}$ & \ct $\geq$1 &$[-+-+]$ or $[+-+-]$ & predator-prey\\
\hline  
\ct $\rm{M}$ &\ct  $\geq$1 &$[++-+]$ & parasitic\\
\hline 
\ct $\rm{C}$ &\ct  $\geq$1 &$[+-++]$ & parasitic\\
\hline 
\ct $\rm{B}_1$ & \ct $\geq$1 &$[++++]$ & cooperative\\
\hline 
\ct $\rm{B}_2$ &\ct  $\geq$1 &$[-+++]$ & cooperative dependency\\
\hline 
\ct $\rm{B}_3$ &\ct  $\geq$1 &$[+++-]$ & cooperative dependency\\
\hline 
\ct $\rm{B}_4$ & \ct $\geq$1 &$[-++-]$ & codependency \\
\hline 
\ct $\rm{RR}$ & \ct $\geq$2 & $[+--+]$ & competitive\\
\hline 
\ct $\rm{GG}$ & \ct $\geq$2 & $[-+-+]$ or $[+-+-]$ & predator-prey\\
\hline 
\ct $\rm{MM}$ & \ct $\geq$2 & $[++-+]$ & parasitic\\
\hline 
\ct $\rm{CC}$ & \ct $\geq$2 & $[+-++]$ & parasitic\\
\hline 
\ct $\rm{RG}$ &\ct  $\geq$2 & $[+--+]$ or $[-+-+]$ or $[+-+-]$ & competitive or~predator-prey\\
\hline 
\ct $\rm{GM}$ & \ct $\geq$2 & $[-+-+]$ or $[+-+-]$ or $[++-+]$ & predator-prey or parasitic \\
\hline 
\ct $\rm{RM}$ & \ct $\geq$2 & $[+--+]$ or $[++-+]$ & competitive or parasitic\\
\hline 
\ct $\rm{GC}$ & \ct $\geq$2 & $[-+-+]$ or $[+-+-]$ or $[+-++]$ & predator-prey or parasitic\\
\hline 
\ct $\rm{RC}$ & \ct $\geq$2 & $[+--+]$ or $[+-++]$ & competitive or parasitic \\
\hline 
\ct $\rm{NE}$ &\ct  0 &  & \\
\hline 
\end{tabular}
\end{table}

As depicted in the $P_2$-diagram (Figure~\ref{figure:P2diagram}), we make the following statements about the solutions of the inverse problem:

\begin{enumerate}
\item Any trajectory of system \eref{LV3} that starts at $P_0$ remains in the first quadrant. Thus, there are no solutions of the inverse problem with $P_2$ outside of the first quadrant, and the $P_2$-diagram is thereby restricted to that quadrant.  
\item The first quadrant can be partitioned into open regions $\mathcal{R}_\Omega$ in which there are solutions to the inverse problem with particular sign structure $\sigma_A$ of $A$. Regions labeled by the same subscript and shown in the same color share the same sign structure for $A$, as indicated in Table \ref{table:P2diagram}. 
\item If $P_2\in \mathcal{R}_{\rm{NE}}=\cup_{j=1}^4 \mathcal{R}_{\rm{NE}_j}$, then the inverse problem has no solution. 
\item  If $P_2$ is located in any labeled region $\mathcal{R}_\Omega$ not included in $\mathcal{R}_{\rm{NE}}$ except for regions $\mathcal{R}_{\rm{G}}$ or regions labeled with two letters, or $P_2$ lies on any boundary between regions represented by a solid curve in the $P_2$-diagram, then the inverse problem has a unique solution.
\item In regions $\mathcal{R}_{\rm{G}}$ the inverse problem has a countable family of solutions that correspond to periodic orbits (discussed below in Section \ref{subsection:non-uniqueness:green}).
\item In regions $\mathcal{R}_\Omega$ labeled by two letters (but not by NE), two solutions 
arise due to a fold in the manifold $\mathcal{M}$.
\item The regions are separated by curves $\mathcal{C}_\omega$ on which some of the parameters vanish (discussed  below in Section \ref{boundarycurves}). 
\end{enumerate}

\begin{figure}
\centering  
\includegraphics[width=1\textwidth]{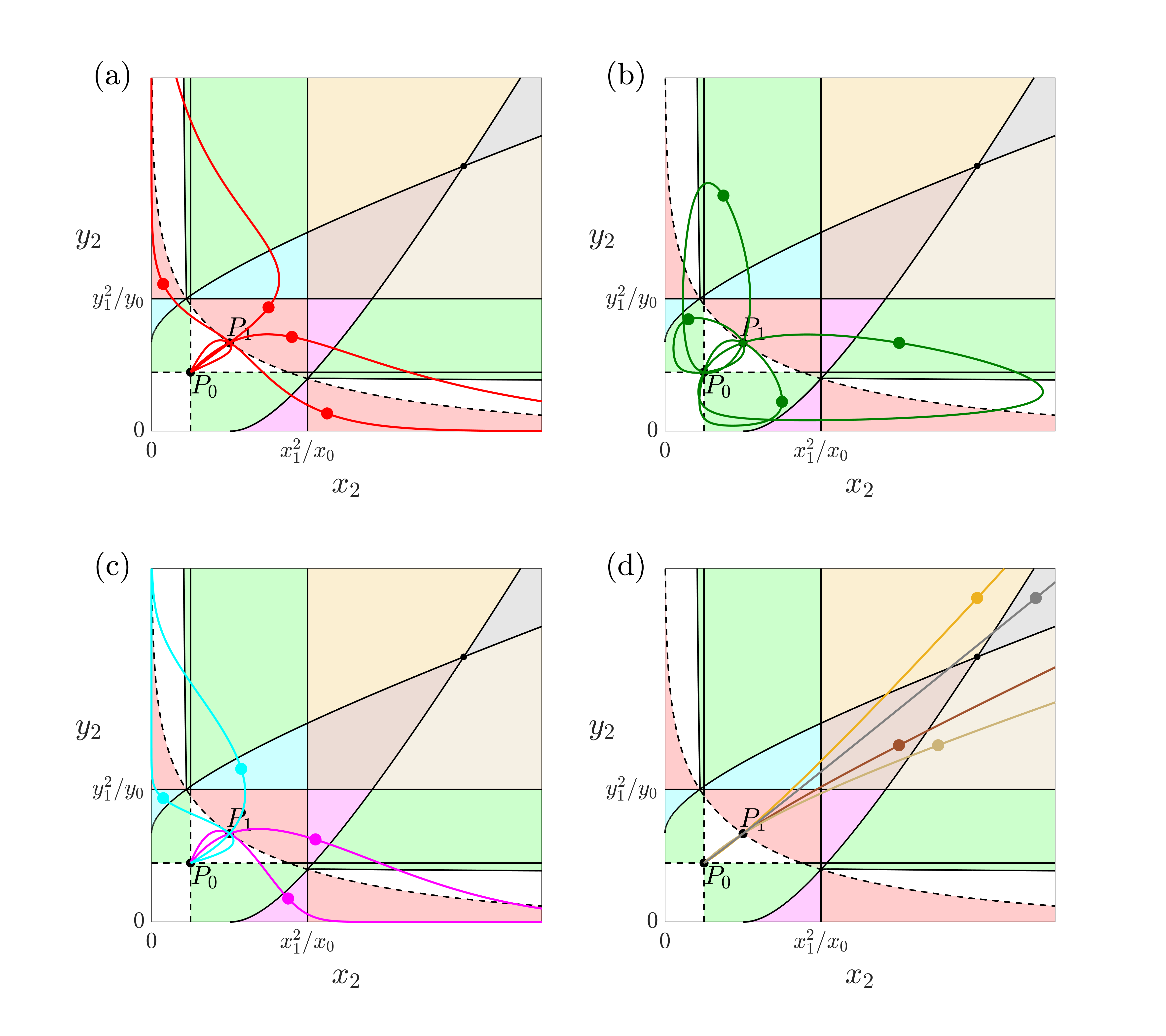}
\caption{Example trajectories with different signatures $\sigma_A$. Colored dots without labels indicate locations of $P_2$ for the illustrated trajectories.}\label{exampletrajs}           
\end{figure}

Figure \ref{exampletrajs}  shows examples of trajectories that pass through specific choices of the point $P_2$ in various existence regions. Notice that each such trajectory passes through several regions in Figure \ref{figure:P2diagram} including the non-existence region, $\mathcal{R}_{\rm{NE}}$.  These behaviors are consistent with the nature of this diagram; that is, the $P_2$-diagram provides information about what we can infer about $A$ from the location of $P_2$ and does not tell us what regions of the $(x,y)$ phase space the trajectory does or does not visit before or after it reaches $P_2$.

\subsection{Trivial non-uniqueness for periodic solutions}
\label{subsection:non-uniqueness:green}
We next go into more detail about the $P_2$-diagram for the LV system (\ref{LV3}), including some aspects that differ from those that can occur for linear and affine systems.
First, we note that due to the presence of a conserved quantity, or Hamiltonian, in the system (see Section 2 of the Supplementary Materials), in some regions of $P_2$ values, there is a form of trivial non-uniqueness of parameter matrices $A$ for which $P_j = \phi(j,A)$ for all $j \in \{ 0, 1, 2\}$.
In the cases when $\sigma_A=[-+-+]$ or $[+-+-]$, each trajectory through $P_0$, $P_1,$ and $P_2$ lies on a periodic trajectory, which is a level set of the Hamiltonian. This gives infinitely many solutions to the inverse problem by the following rescaling argument: Let $A$ be a parameter matrix that solves the inverse problem for data $d=\{P_0,P_1,P_2\}$, with $\sigma_A=[-+-+]$ or $[+-+-]$ and trajectory $\phi(t)$.  Let $T$ be the period of the orbit and $\tau$ be the time of the shortest clockwise passage from $P_0$ to $P_1$ (or shortest counterclockwise passage from $P_1$ to $P_0$). If the trajectory with $\phi(0)=P_0$ travels $m$ times clockwise around the orbit before $\phi(1)=P_1$, then we have $1=mT+\tau$. The same relation holds for a trajectory with $\phi(0)=P_0$ that travels counterclockwise, with $m=-1$ for a trajectory that reaches $\phi(1)=P_1$ before completing a full orbit, and $m<-1$ for a trajectory that travels counterclockwise $1-m$ full orbits before reaching $\phi(1)=P_1$. Autonomy of system \eref{LV3} implies that the same relation also holds for travel from $P_1=\phi(1)$ to $P_2=\phi(2)$. Consider now a system \eref{LV3} (i.e., \eref{LIP}) with the parameter matrix $\tilde{A}=\gamma A$. Since the system is linear-in-parameters, the orbit of the trajectory starting at $P_0$ will not change but the period of the trajectory will change to $\tilde{T}=T/\gamma$ and the shortest time between $P_0$ and $P_1$ to $\tilde{\tau}=\tau/\gamma$. If $\gamma = (\tilde{m}-m)T+1$ for some integer $\tilde{m}$ then the trajectory of the system with parameter matrix $\tilde{A}$ again obeys $\phi(i)=P_i, i=0,1,2$ and hence $\tilde{A}$ is another solution of the inverse problem, one which travels $\tilde{m}$ times clockwise around the orbit. 
This observation explains why in Table \ref{table:P2diagram} the cases $[-+-+]$ and $[+-+-]$ always show up in pairs. 

We treat this type of non-uniqueness as trivial and count all of these matrices as a single solution to the inverse problem, which we represent with the parameter matrix $A$ for which $m=0$ (clockwise trajectory) or $m=-1$ (counterclockwise trajectory), depending on which of those does not travel a full orbit before reaching $P_2$ when starting from $P_0$.
Moreover, as the periodic orbit is convex
(see Section 2 of Supplementary Materials) this choice uniquely determines the signature of the inverse problem solution to be of predator-prey type, i.e., we have $\sigma_A=[- + - +]$ (and the trajectory travels clockwise) if $\phi(2)$ is below the straight line $P_0P_1$, and $\sigma_A=[+ - + -]$ (and the trajectory travels counterclockwise)  if $\phi(2)$ is above that line.
Note, however, that in contrast to this trivial form of non-uniqueness, the  $\mathcal{R}_{\rm{GG}}$ region in the $P_2$-diagram and in Table \ref{table:P2diagram} corresponds to more than one solution with $\sigma_A=[-+-+]$ or $[+-+-]$; here, the multiple trajectories truly represent different orbits.
The trivial form of non-uniqueness and the argument we used to establish its presence carry over to periodic solutions of linear and affine systems, since we can use the same multiplicative parameter rescaling idea in those settings as well.  
The non-trivial non-uniqueness of $\mathcal{R}_{\rm{GG}}$, on the other hand, is a new feature not found in linear and affine systems.



\subsection{Boundary curves in the $P_2$-diagram}\label{boundarycurves}
\begin{figure}
\centering  
\includegraphics[width=1\textwidth]{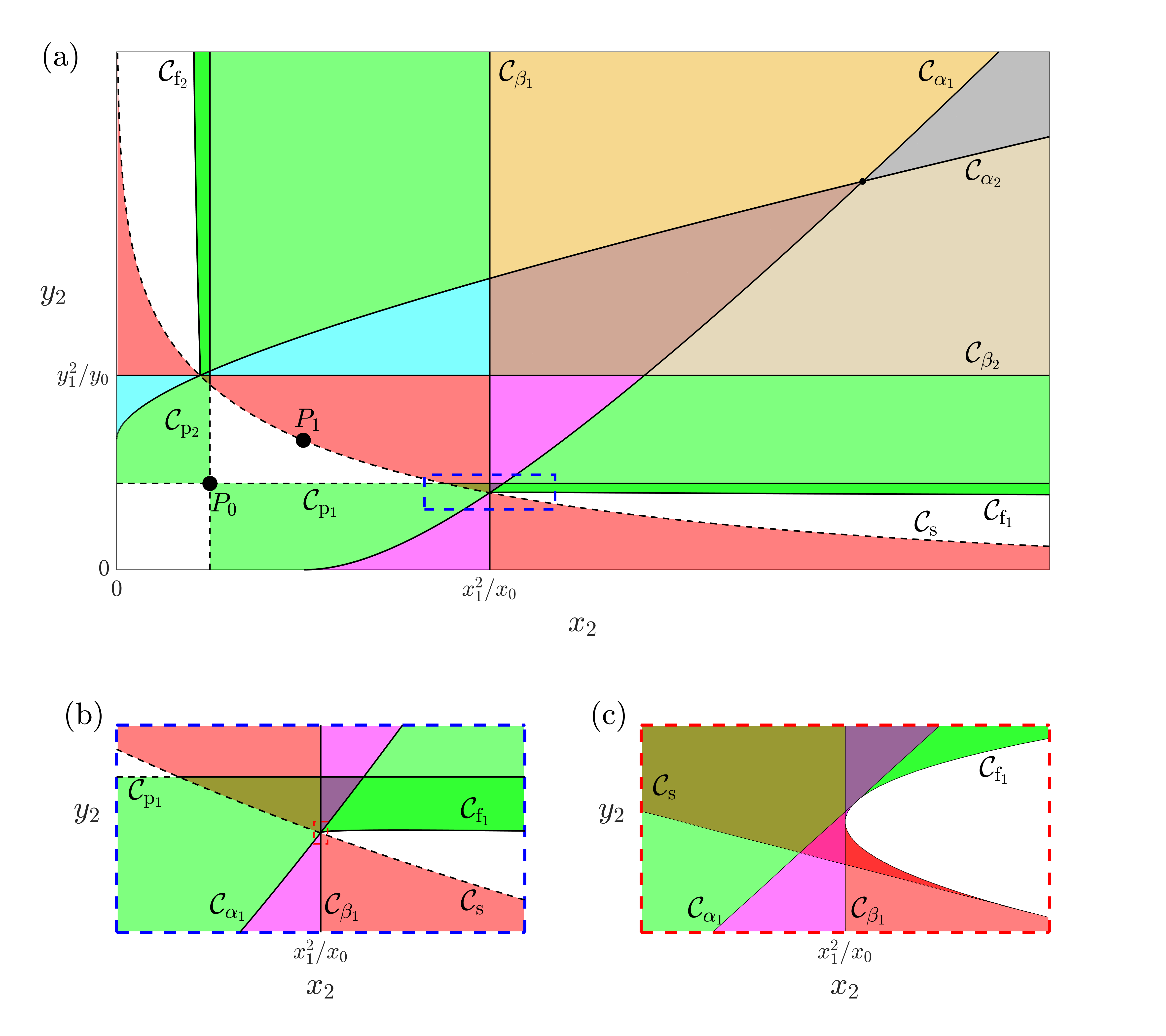}
\caption{The $P_2$-diagram, which is also displayed in Figure \ref{figure:P2diagram}, showing the labels of the curves introduced in this section.}
\label{figure:P2diagram:curvenames}           
\end{figure}

A key feature of the $P_2$-diagram is that the regions with different sign signatures or identifiability properties are separated by smooth, simple boundary curves.  Thus, a full specification of the $P_2$-diagram requires an explanation of what these boundary curves represent, in terms of model parameters.  Therefore, we next focus on these curves, which we label in Figure~\ref{figure:P2diagram:curvenames} and  classify into three different types: (i) the curves along which at least one of the parameters $\aq,\bq,\bw,\aw$ is zero and which therefore separate the regions of the diagram with different $\sigma_A$; (ii) the curves along which at least one of the parameters $\aq,\bq,\bw,\aw$ approaches infinity; and (iii) the curves that describe folds in the manifold $\mathcal{M}$ resulting in different solutions $A$ that co-exist for the same value of $P_2$.  Most of the curves are found numerically and in Section 1 of Supplementary Materials
we discuss the methods used in this section.
We note that the overlap of regions in $P_2$-space in which solutions with different dynamical behavior co-exist is a new feature that arises in this LIP setting, relative to the linear and affine cases, and thus represents an important example of the complexity that nonlinearity can add to solving the inverse problem.

\subsubsection{Signature boundaries $\mathcal{C}_{\bq}$, $\mathcal{C}_{\bw}$, $\mathcal{C}_{\aq}$, $\mathcal{C}_{\aw}$}\label{subsection:Curve:parais0}
The boundaries separating regions of fixed $\sigma_A$ in the $P_2$-diagram correspond to the cases when at least one of the parameters $\aq,\bq,\bw,\aw$ is zero. For certain combinations of zero parameters one can find an analytical solution of system \eref{LV3}. The detailed derivations of solutions for these cases are given in Section 3 of Supplementary Materials.
We summarize those results below.

When $\bq=0$, the $x$-equation in system \eref{LV3} becomes decoupled and has the trajectory $x(t)=\rme^{\aq t}x_0$.  Hence, $\aq=\ln(x_1/x_0)$ and $x_2=x_1^2/x_0$, so $P_2$ for such a system lies on a vertical straight line $\mathcal{C}_{\bq}$ in the $(x_2,y_2)-$plane (see Figure \ref{figure:P2diagram:curvenames}). We can determine explicitly the parameter values corresponding to any point on $\mathcal{C}_{\bq}$ as functions of $P_0,P_1,P_2$. For example, the intersection between $\mathcal{C}_{\bq}$ and $\mathcal{C}_{\aw}$ is given by
\begin{equation}
    \label{eq:beta1=0}
P_{\bq,\aw}=\left(x_1\left(\frac{x_1}{x_0}\right), y_1\left(\frac{y_1}{y_0}\right)^\frac{x_1}{x_0}\right).
\end{equation}

An analogous situation holds in the case when $\bw=0$, with $P_2$ then lying on the horizontal line $y=y_1^2/y_0$, which we denote $\mathcal{C}_{\bw}$. The intersection between $\mathcal{C}_{\bw}$ and $\mathcal{C}_{\aq}$ is given by
\begin{equation}
\label{eq:beta2=0}
P_{\bw,\aq}=\left(x_1\left(\frac{x_1}{x_0}\right)^\frac{y_1}{y_0}, y_1\left(\frac{y_1}{y_0}\right)\right).
\end{equation}

We have not been able to find explicit formulas for the curves $\mathcal{C}_{\aq}$ and $\mathcal{C}_{\aw}$.
When $\aq=0$, the first equation of the system becomes $\dot{x}=\bq xy$, and $x_1>x_0$ implies $\beta_1>0$.  We use numerical computation (see Section 1 of Supplementary Materials) 
to obtain the curve $\mathcal{C}_{\aq}$ of points $P_2$ corresponding to $\aq=0$. The curve $\mathcal{C}_{\aq}$ appears to be single-valued in $y_2$ and defined for $y_2\in (0,\infty)$.
Since $\bq>0$, this curve lies to the right of $x_1$.   Similarly, we can find numerically the curve $\mathcal{C}_{\aw}$ corresponding to $\aw=0$.

When both $\aq=\aw=0$, however, the system \eref{LV3} can be rewritten as a Bernoulli equation for $x$, and the trajectory is a straight line passing through $P_0$ and $P_1$. The intersection point of the curves $\mathcal{C}_{\aq}$ and $\mathcal{C}_{\aw}$ is
\[
 P_{\aq,\aw}=\left(x_1\left(\frac{x_1}{x_0}\right),
 y_1\left(\frac{y_1}{y_0}\right)\right)\left(\frac{x_1}{x_0}+\frac{y_1}{y_0}-\frac{x_1y_1}{x_0y_0}\right)^{-1} .
\]
Note that the scaling factor multiplying the expression is of the form $(a+b-ab)^{-1}$ where $a=\frac{x_1}{x_0}>1$ and $b=\frac{y_1}{y_0}>1$ for our choice of the data points $P_0,P_1$. 
The scaling factor is either larger than one or negative under these conditions.  In the former case, by comparison with the expressions \eref{eq:beta1=0}, \eref{eq:beta2=0}, we observe that  $P_{\aq,\aw}$ is to the right of the line $\mathcal{C}_{\bq}$ and above $\mathcal{C}_{\bw}$ as in Figures~\ref{figure:P2diagram} and \ref{figure:P2diagram:curvenames}. On the other hand, as $a$ and $b$ grow, the scaling factor and hence the intersection point blow up to infinity, beyond which  the ratios are large so that $a^{-1}+b^{-1}<1$,  the scaling factor is negative, and the curves $\mathcal{C}_{\aq}$ and $\mathcal{C}_{\aw}$ do not intersect in the positive quadrant. One implication of these observations is that the $P_2$-diagram is not invariant under a change in the location of the points $P_0,P_1$ even if that change preserves the inequalities between their $x$ and $y$ coordinates. 

Note that in the $P_2$-diagram, an intersection occurs between 
$\mathcal{C}_{\aq}$ and $\mathcal{C}_{\bq}$. However, any solution of the system \eref{LV3} with $\aq=0$ and $\bq=0$ has $x(t)=x_0$ and hence cannot pass through both of the specified points $P_0$ and $P_1$.  Thus, the apparent intersection point of these curves is only due to the projection of $\mathcal{M}$ onto the $P_2$-plane. Such an intersection necessarily corresponds to two distinct solutions and hence must be a point of non-uniqueness for solutions to the inverse problem. We discuss more details about non-uniqueness below.


\subsubsection{Separatrix $\mathcal{C}_{\rm{s}}$}

The separatrix $\mathcal{C}_{\rm{s}}$ is the part of the boundary of the red region defined by a set of locations in the $P_2$ plane at which the solution of the inverse problem ceases to exist. Numerical computations show that as $\mathcal{C}_{\rm{s}}$ is approached from the direction of  $\mathcal{R}_{\rm{R}}$, the solutions of the inverse problem are such that all four parameters $\aq,\bq,\bw,\aw$ diverge, with $\alpha_1, \alpha_2 \to \infty$ and $\beta_1, \beta_2 \to -\infty$. Indeed, as $P_2$ approaches $P_1$ through $\mathcal{R}_{\rm{R}}$, the parameters must scale in such a way that the overall rates of change of $x$ and $y$ near $P_1$ become arbitrarily small, corresponding to the shrinking distance between $P_1$ and $P_2$ that is covered in one time unit. It is difficult to compute $\mathcal{C}_{\rm{s}}$ numerically because one cannot simply proceed with continuation along a specific parameter.

The separatrix $\mathcal{C}_{\rm{s}}$ has an apparent intersection with $\mathcal{C}_{\bq}$ that lies below the horizontal line $y=y_0$. Similarly, $\mathcal{C}_{\rm{s}}$ seems to intersect the line $\mathcal{C}_{\bw}$ to the left of $x=x_0$.
Intersections between a curve where parameter values diverge and curves where individual parameter values are 0 suggest the existence of some non-trivial three-dimensional structure to $\mathcal{M}$, which we discuss further below. The graph of $\mathcal{C}_{\rm{s}}$ is shown in panel (a) of Figure \ref{figure:P2diagram:curvenames} as the dashed line going from the top left corner to the bottom right corner, forming the boundary between the red regions $\mathcal{R}_{\rm{R}}$ and the white regions $\mathcal{R}_{\rm{NE}_4},\mathcal{R}_{\rm{NE}_2},$ and $\mathcal{R}_{\rm{NE}_3}$ (see Figure \ref{figure:P2diagram} for region labels).

\subsubsection{Periodic orbit limits $\mathcal{C}_{\rm{p}_1}$, $\mathcal{C}_{\rm{p}_2}$}

\begin{figure}
\centering 
\includegraphics[width=0.9\textwidth]{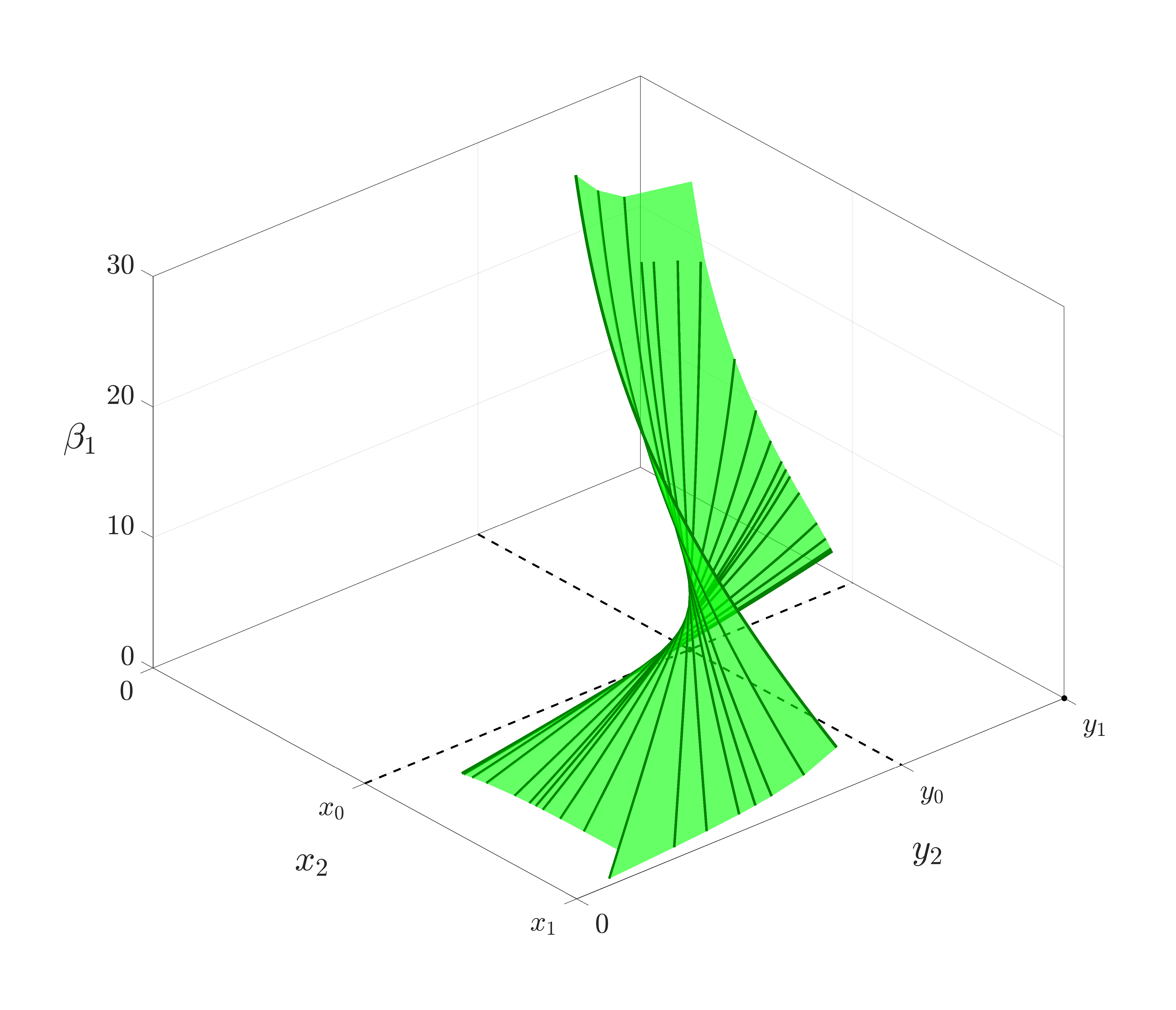}
\caption{ 
The structure $\mathcal{M}_{\bq}$ formed by using the $\bq$ values along  $\mathcal{M}$ to embed its projection over the region $(0,x_1)\times (0,y_1)$ into $(x_2,y_2,\bq)$-space. The dark green curves  correspond to lines passing through $P_0$ with different $x/y$ slopes and $\bq$ heights. $\mathcal{M}_{\bq}$ behaves like a ruled surface near $P_0$. Note the divergence of $\bq$ as $x_2$ approaches $x_0$ and as $y_2$ approaches $y_0$. }\label{figure:3Db1vsP2_P2nearP0}           
\end{figure}

There are two curves in the $P_2$-diagram labeled $\mathcal{C}_{\rm{p}_1}$, $\mathcal{C}_{\rm{p}_2}$, which form  boundaries of the green regions of the $P_2$-diagram 
that correspond to periodic orbits. Several parameters go to infinity when $P_2$ approaches 
the curve $\mathcal{C}_{\rm{p}_1}$ from below or the curve $\mathcal{C}_{\rm{p}_2}$ from the left. They are discussed below and shown in Figures \ref{figure:3Db1vsP2_P2nearP0}, \ref{figure:Fold},  and \ref{figure:3Db1vsP2_ver1}. For example,  in Figure \ref{figure:3Db1vsP2_P2nearP0}, we take $\bq$ as a representative parameter and plot a two-dimensional structure $\mathcal{M}_{\bq}$ that corresponds to the $P_2$-diagram near $P_0$ with $\bq$ as a third coordinate, which illustrates the divergence of $\bq$ along certain curves with $\{y=y_0\}$.  


\subsubsection{The fold curves $\mathcal{C}_{\rm{f}_1}$, $\mathcal{C}_{\rm{f}_2}$}
We have found that there are two folds 
in the projection of the manifold $\mathcal{M}$  onto the $P_2$-diagram, and each of these folds forms a boundary between a region of $P_2$ on which the inverse problem has two solutions in parameter space and a region with no inverse. To obtain a more detailed picture of each fold we fix $x_2>x_1^2/x_0$ and numerically compute the value of $\bq$ as a function of $y_2$ using numerical methods described in Section 1 of the Supplementary Materials.
In Figure \ref{figure:Fold}, a point on the fold curve is represented by a blue dot with $y_2 = y_F < y_0$. On the bottom branch of the diagram, $\bq$ decreases as $y_2$ increases, while on the top branch, $\bq$ increases to $\infty$ as $y_2$ approaches $y_0$ (and $(x_2,y_2)$ approaches $\mathcal{C}_{\rm{p}_1}$). The fold results in non-uniqueness of the solution of the inverse problem for $y_2 \in (y_F, y_0).$ Both branches of solutions have $\sigma_A=[-+-+].$

\begin{figure}
\centering  
\includegraphics[width=1\textwidth]{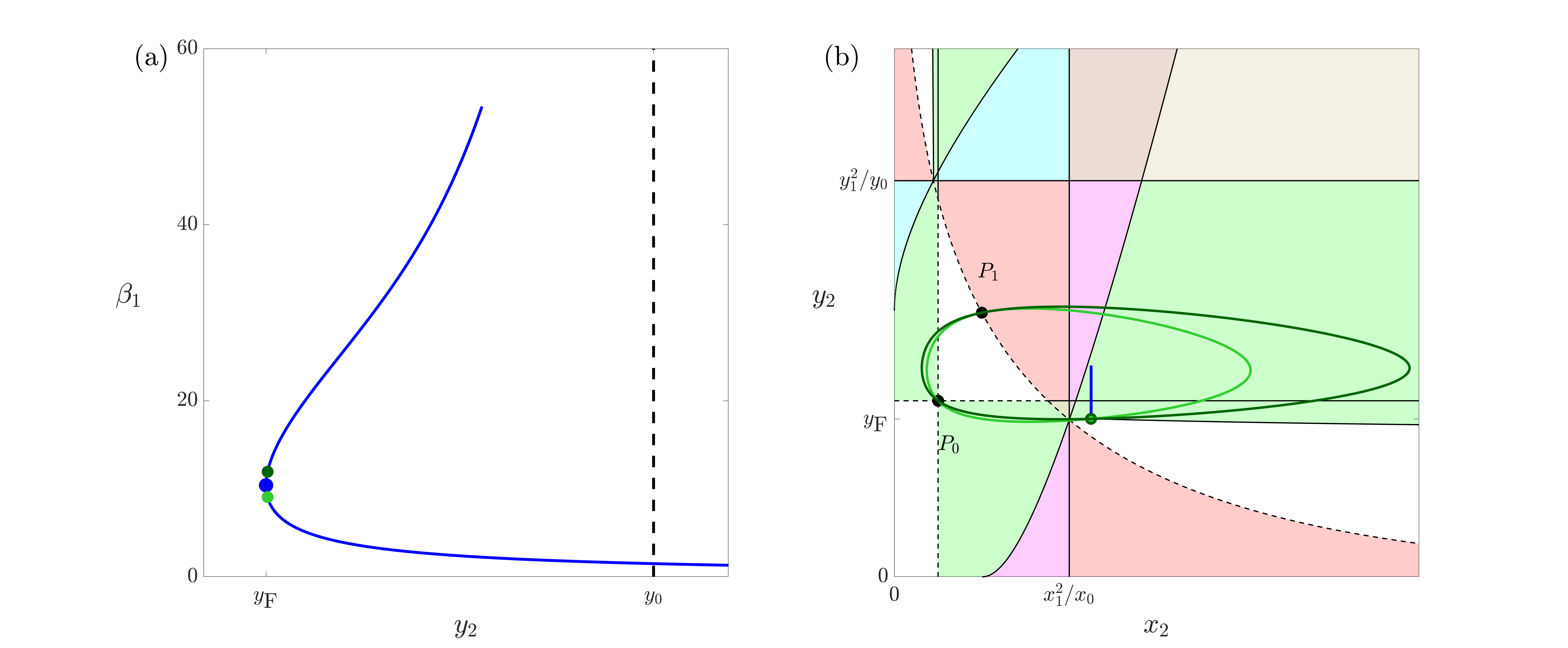}
\caption{ (a) Cross-section of the manifold $\Mbq$ at fixed $x_2 = 4.5$ (with $P_0=[1,1]$ and $P_1=[2,1.5]$) near a fold, located at $(\bq,x_2,y_2) = (\bq,x_2,y_F)$. The blue dot is the fold, while the deep and light green dots are nearby points on $\Mbq$ that share the same $y_2$ coordinate, slightly greater than $y_{\rm{F}}$. (b) The $P_2$-diagram augmented with the two trajectories of system \eref{LV3} that solve the inverse problem when $x_2=4.5$ and $y_2$ takes the value indicated by the green dots in (a). These dots lie on top of each other in (b).  The blue ray from the green dot is a segment in $\{x_2=4.5\}$ corresponding to the range of $y_2$ values, from just below $y_{\rm{F}}$ to above $y_0$, for which the cross-section in (a) is displayed. }
\label{figure:Fold}     
\end{figure}

The curve $\mathcal{C}_{\rm{f}_1}$  obtained as a two-parameter continuation in the
$(x_2,y_2)$ diagram of the fold point shown in Figure \ref{figure:Fold}(a) is one of the fold curves
and appears in the  $P_2$-diagram in Figure 
\ref{figure:P2diagram:curvenames}.  On the manifold $\mathcal{M}$ the curve 
$\mathcal{C}_{\rm{f}_1}$ passes through the regions $\mathcal{R}_{\rm{G}}$, 
$\mathcal{R}_{\rm{M}}$, and $\mathcal{R}_{\rm{R}}$. When projected onto the $P_2$-diagram 
those regions overlap and each overlap gives rise to two possible solutions of the 
inverse problem. As noted previously, in the $P_2-$diagram we denote the overlaps as 
$\mathcal{R}_{\Omega_1\Omega_2}$, where each $\Omega_i$ is one of the letters 
$\rm{R},\rm{G},\rm{M},$ and $\rm{C}$ as in Table \ref{table:P2diagram}. Two distinct 
letters imply that in that overlap region the signature $\sigma_A$ for the inverse 
problem can attain two distinct values, i.e., there are two distinct types of systems compatible with points $P_2$ located in that region.

Notice that there are two fold curves shown in Figure \ref{figure:P2diagram:curvenames}, $\mathcal{C}_{\rm{f}_1}$ and $\mathcal{C}_{\rm{f}_2}$. Fold curve $\mathcal{C}_{\rm{f}_1}$ borders the region $\mathcal{R}_{\rm{NE}_3}$ and the regions labeled $\mathcal{R}_{\rm{R}},\mathcal{R}_{\rm{G}},$ and $\mathcal{R}_{\rm{M}}$, which is sketched in Figure \ref{figure:P2diagram:curvenames} 
panels (a), (b) and (c). Fold curve $\mathcal{C}_{\rm{f}_2}$ is in a symmetric position in the top left part of the diagram beside $\mathcal{R}_{\rm{NE}_4}$ and the regions labeled $\mathcal{R}_{\rm{R}},\mathcal{R}_{\rm{G}},$ and $\mathcal{R}_{\rm{C}}$; the latter is only sketched in panel (a). As is evident in panel (c) of Figure \ref{figure:P2diagram:curvenames}, the projection of the fold curve $\mathcal{C}_{\rm{f}_1}$ onto the $P_2$ diagram intersects tangentially with the projections of the curves $\mathcal{C}_{\aq}$, $\mathcal{C}_{\bq}$, and $\mathcal{C}_{\rm{s}}$.


\subsection{The visualization of the $P_2$-diagram in three dimensions}
\label{sec:vis}
\begin{figure}
\centering  
\includegraphics[width=1\textwidth]{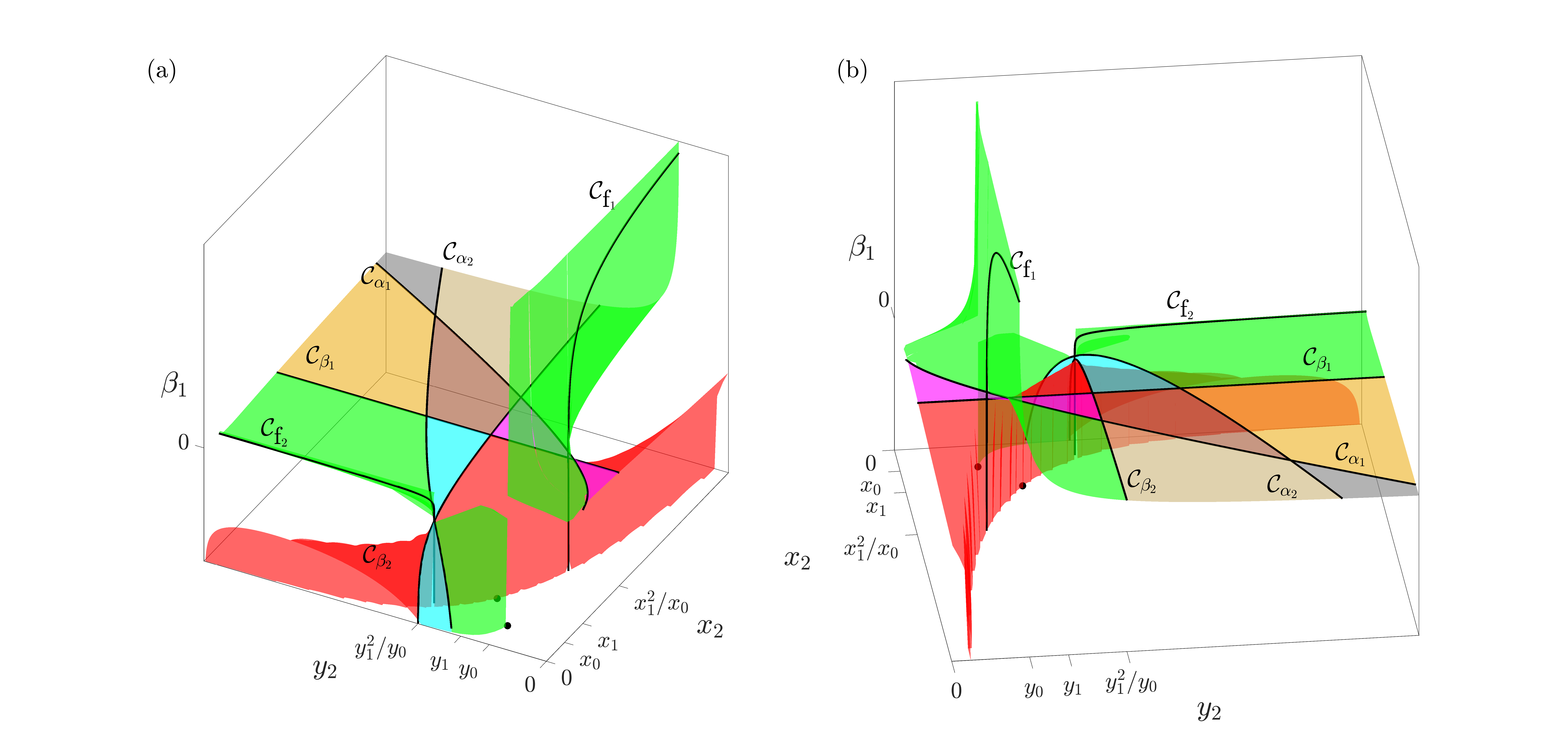}
\caption{\small Two different three-dimensional visualizations of $\mathcal{M}_{\bq}$. These views highlight the intricate folding of the surfaces of inverse problem solutions that give rise to trajectories passing with appropriate timing through the fixed data points $P_0, P_1$ and the additional data point $P_2$ represented in the diagram, with the usual color coding based on $\sigma_A$.
See Section \ref{sec:vis} for more details.}
\label{figure:3Db1vsP2_ver1}           
\end{figure}

\begin{figure}
\centering  
\includegraphics[width=1\textwidth]{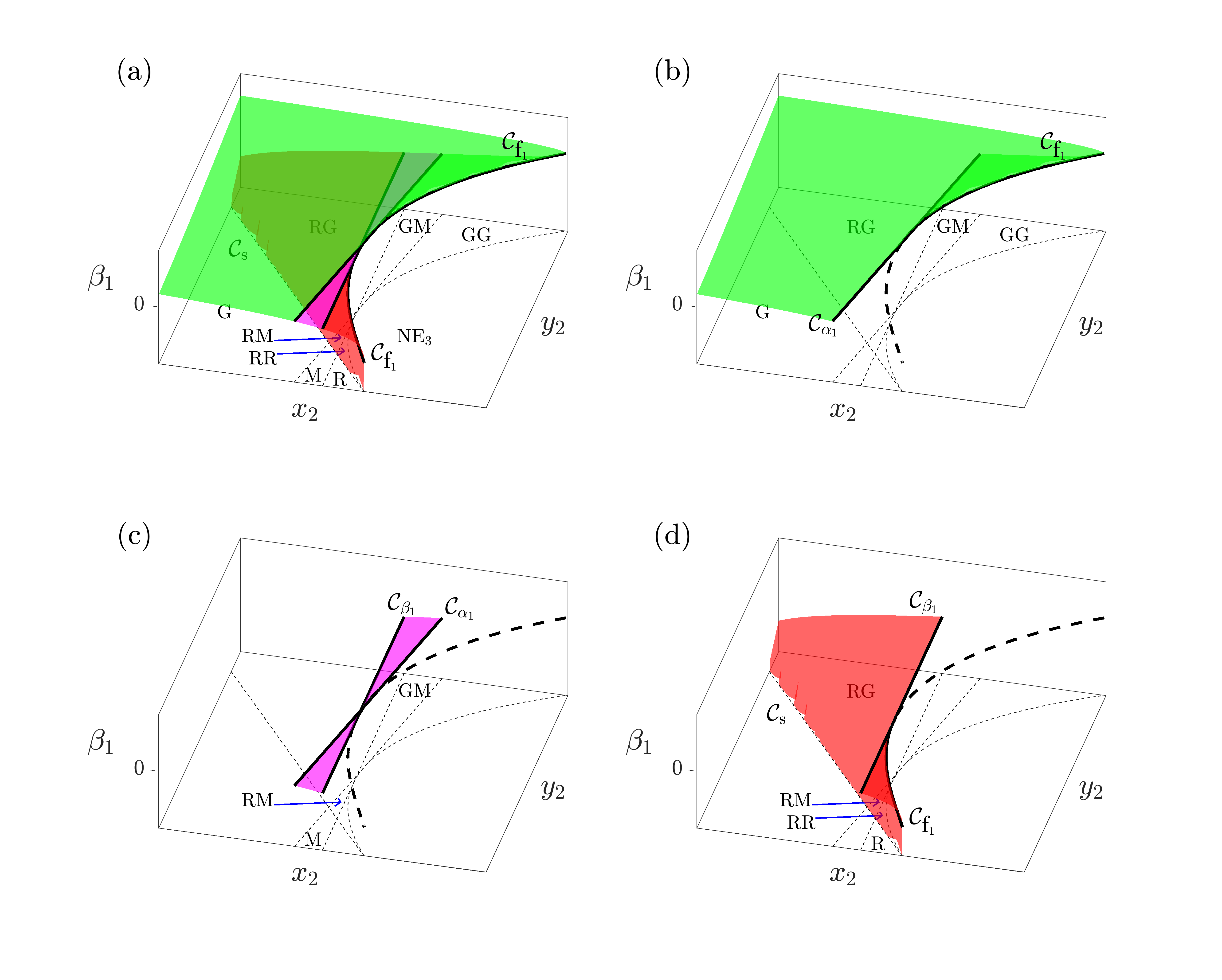}
\caption{Zoomed view of Figure \ref{figure:3Db1vsP2_ver1} focusing on the region in the $P_2$-diagram, as viewed in $(x_2,y_2,\bq)$-space, where the four curves $\mathcal{C}_{\aq}$, $\mathcal{C}_{\bq}$, $\mathcal{C}_{\rm{f}_1}$, and $\mathcal{C}_{\rm{s}}$ interact. The curve $\mathcal{C}_{\rm{f}_1}$ is sketched solid when it participates in the fold in the corresponding panel(s), and dashed when it does not. Other than labeling these curves, we do not include labels on the plotted surface to avoid clutter. We do, however, label the regions in $(x_2,y_2)$-space along the base of each panel according to the color code for the sign signatures that arise there (although we omit the MM label due to the small size of its region; see Figure \ref{figure:P2diagram}). (a) The full collection of color-coded surface components in the zoomed region.  (b)-(d) Separate plots of the (b) green, (c) magenta, and (d) red components of the surface.}
\label{figure:3Db1vsP2_neartwist}           
\end{figure}

Now that we have captured the fold curves along $\mathcal{M}$, we continue to  provide a more complete illustration  of the folds and twists in the $P_2$-diagram, as shown in Figures \ref{figure:3Db1vsP2_P2nearP0}, \ref{figure:3Db1vsP2_ver1},  and \ref{figure:3Db1vsP2_neartwist}. 
To this end,  we take $\bq$ as a representative parameter and plot the surface $\mathcal{M}_{\bq}$ of solutions of the inverse problem in $(x_2,y_2,\bq)$-space.

Figure \ref{figure:3Db1vsP2_ver1} presents two different views of $\mathcal{M}_{\bq}$.  We show the images in $\mathcal{M}_{\bq}$ of some of the curves we described in the previous subsection, namely $\mathcal{C}_{\aq}$, $\mathcal{C}_{\bq}$, $\mathcal{C}_{\bw}$, $\mathcal{C}_{\aw}$, $\mathcal{C}_{\rm{f}_1}$, and $\mathcal{C}_{\rm{f}_2}$. The separatrix $\mathcal{C}_{\rm{s}}$ is not visible as a curve in this figure since it corresponds to $\bq\to-\infty.$ We see that when $P_2\in[x_1^2/x_0,\infty)\times[y_1^2/y_0,\infty)$, the surface is relatively flat. The surface twists over near $P_0$; we have already seen in Figure \ref{figure:3Db1vsP2_P2nearP0}  that when $P_2$ is exactly at $P_0$, there is non-uniqueness of the solutions of the inverse problem and the distribution of the solutions is no longer discrete; we will discuss this phenomenon in the next subsection.  The surface also folds over in the region corresponding to panel (c) of Figures \ref{figure:P2diagram} and \ref{figure:P2diagram:curvenames} (that is, the region where the four curves
$\mathcal{C}_{\aq}$, $\mathcal{C}_{\bq}$, $\mathcal{C}_{\rm{f}_1}$, and $\mathcal{C}_{\rm{s}}$ interact), so we sketch zoomed views focusing on  these two regions separately in Figures \ref{figure:3Db1vsP2_P2nearP0} and  \ref{figure:3Db1vsP2_neartwist}, respectively.  In these regions, when we color parts of $\mathcal{M}_{\bq}$ according to their sign signatures $\sigma_A$, we in fact have two surfaces to consider.  Within each surface, we apply our usual convention of using red for $\sigma_A=[+--+]$, magenta for $\sigma_A=[++-+]$, and cyan for $\sigma_A=[+-++]$, as can be seen in Figures \ref{figure:3Db1vsP2_ver1},  \ref{figure:3Db1vsP2_neartwist}. In Figure \ref{figure:3Db1vsP2_neartwist} specifically, within the $(x_2,y_2)$-plane that forms the base of each panel, we label each of several regions with the pair of letters indicating those $\sigma_A$ that arise on these regions.

\subsection{Non-uniqueness when $P_2=P_0$}

From the previous subsections, we see that each fold curve forms the boundary of a region of non-uniqueness of the solution to the inverse problem.
Recall from Section \ref{subsection:non-uniqueness:green} that there is another source of non-uniqueness that we treat as trivial, corresponding to simple time rescaling, when $\sigma_A=[-+-+]$ or $[+-+-]$. 

Now we study a special case of non-uniqueness when $P_2=P_0$. This case also yields $\sigma_A=[-+-+]$ or $[+-+-]$, but as we have seen in Figure \ref{figure:3Db1vsP2_P2nearP0}, here there exist infinitely many non-trivial solutions of the inverse problem, which form a continuum in parameter space; that is, in Figure \ref{figure:3Db1vsP2_P2nearP0}, the surface  $\mathcal{M}_{\bq}$
contains a vertical line over $P_0,$
where $\{x_2 = x_0 \}$ intersects $\{ y_2 = y_0 \}$.

We can explain this special non-uniqueness from the following perspective.
There is a family  $\mathcal{T}$ of trajectories of system \eref{LV3} starting at $P_0$ and passing through $P_1$ at $t=1$, but only one of them contains a trajectory passing through $P_2$ at $t=2.$   When $P_2=P_0$, only two distinct points on the trajectory are specified, and hence the cardinality of the family $\mathcal{T}$ is larger, so that the number of trajectories passing through $P_2(=P_0)$ at $t=2$ becomes infinite. 


Notice that in this case, the period of each of the elliptic orbits is two, which means that starting at $P_0$, it takes 1 time unit for the flow to evolve to $P_1$ and an additional 1 time unit for the flow to evolve from $P_1$ to $P_2$.  Thus, if we have a solution $A$, then the matrix $-A$ (trivially) provides another solution, as discussed in Section \ref{subsection:non-uniqueness:green}.

\subsection{$P_2$-diagrams for alternative positions of $P_0$ and $P_1$}\label{sec:altP2}

\begin{figure}
\centering  
\includegraphics[width=1\textwidth]{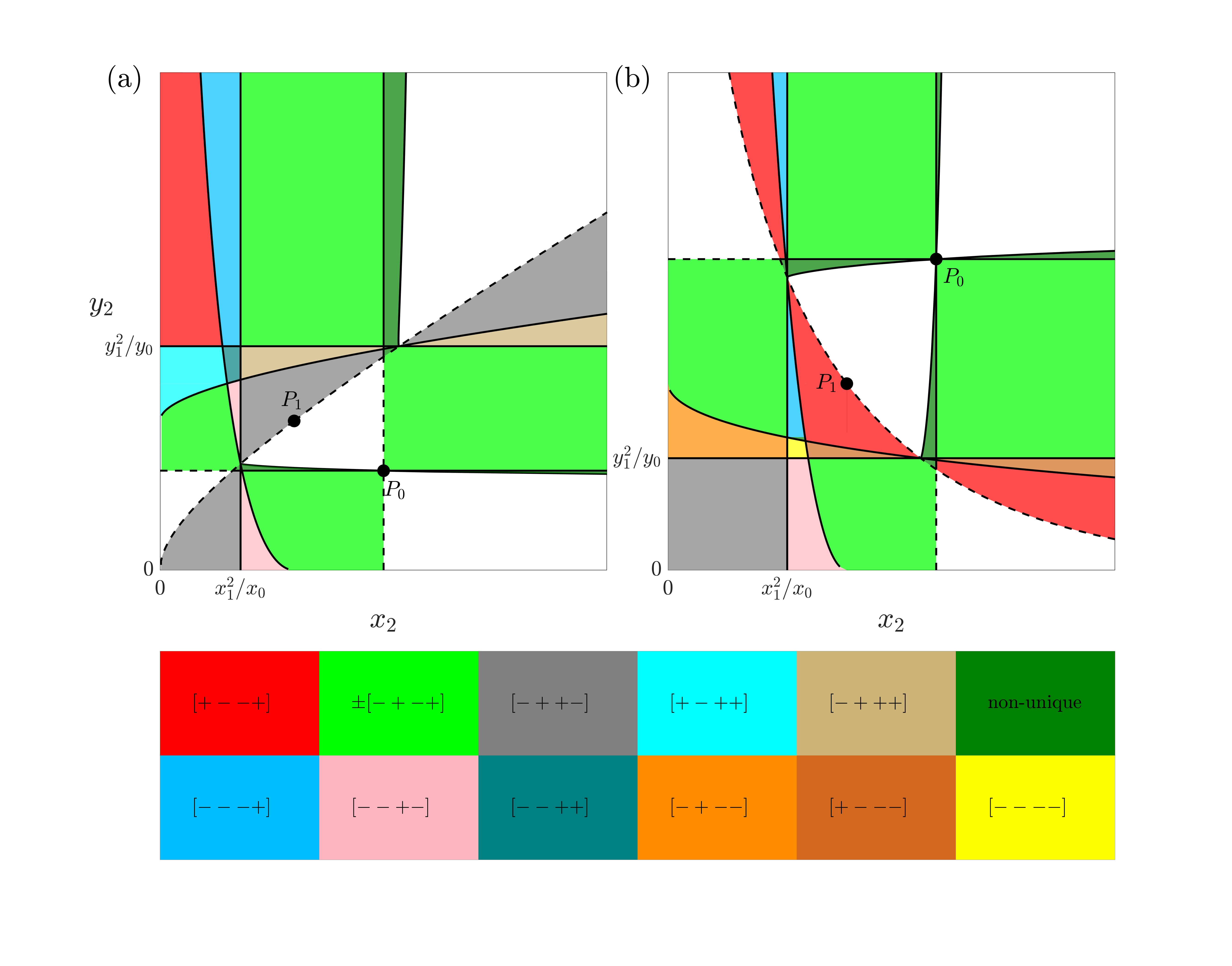}
\caption{The $P_2$-diagrams  in the cases where (a) $x_1<x_0$ and $y_1>y_0$ and (b) $x_1<x_0$ and $y_1<y_0$. The $\sigma_A$ for each region is as labeled in the color chart below the main plots. }\label{figure:DiagramTLBL}
\end{figure}

In previous subsections we focused in detail on the situation where $P_0$ and $P_1$ are positioned so that condition $(\rm{C}_1)$ applies: $x_0<x_1$ and $y_0<y_1$. For completeness and comparison, Figure \ref{figure:DiagramTLBL} shows the $P_2$ diagram for cases in which the positions of $P_0$ and $P_1$ obey alternative conditions $(\rm{C}_2)$: $x_1<x_0$ , $y_1>y_0$, or $(\rm{C}_3)$: $x_1<x_0$, $y_1<y_0$. With the new conditions, new possible signatures $\sigma_A$ appear for solutions of the inverse problem, and different regions now correspond to the signatures we have already discussed above, but there are also many similarities between the new diagrams and the case when $(\rm{C}_1)$ holds.  
For example, in these $P_2$-diagrams, the curves $\mathcal{C}_{\alpha_i}$, $\mathcal{C}_{\beta_i}$, folds $\mathcal{C}_{\rm{f}}$, and separatrices $\mathcal{C}_{\rm{s}}$ also exist, and the lines $\mathcal{C}_{\beta_i}$ are located at the same coordinate values, i.e., $x=x_1^2/x_0$ and $y=y_1^2/y_0$. The nonexistence region in both cases again separates into four disjoint components. There are again regions of non-trivial non-uniqueness, which for simplicity we color dark green in Figure \ref{figure:DiagramTLBL}, irrespective of which types of solutions co-exist.  We do not show a diagram for the case when $x_1>x_0$ and $y_1<y_0$ as it is equivalent to  $(\rm{C}_2)$ with $x$ and $y$ interchanged.

\section{Non-conservative Lotka-Volterra systems}
\label{sec:noncon}

The Lotka-Volterra system \eref{LV3} that we analyzed in previous sections is a conservative system; that is, there is a quantity, the Hamiltonian, that is conserved along trajectories. As a result, we have a trivial non-uniqueness of solutions of the inverse problem corresponding to periodic trajectories, where there are countable families of trajectories that travel along identical closed orbits, as discussed in Section \ref{subsection:non-uniqueness:green}. Linear and affine systems, discussed in Sections \ref{linearresults} and \ref{affineresults}, do not have such trivial non-uniqueness and instead non-unique solutions of the inverse problem give rise to spiraling solutions with various numbers of rotations occurring between the times when the trajectory passes through the data points. We show in the following sections that removal of the conservative nature of the system eliminates the trivial non-uniqueness observed for \eref{LV3}.

\subsection{Lotka-Volterra system with rotated field}
One way in which we can take away the conservative nature of the system is by introducing a rotation-like perturbation that mixed the terms in the vector field \cite{duff1953limit,perko1993rotated}.  Specifically, we consider the system  
\begin{equation}             
\left\{
             \begin{array}{lr}
             \dot{x}=\aq x + \bq x y - p ( \aw y + \bw x y), &  \\
             \dot{y}=\aw y + \bw x y + p ( \aq x + \bq x y), &  \\
             \end{array}
\right.
\label{LVrotate}
\end{equation}
where $p$ is a small parameter, which reduces to system \eref{LV3} when $p=0$, and which vector field is rotated at each point by a fixed angle with respect to the field of \eref{LV3}. As a result, the trajectories of \eref{LVrotate} intersect transversally the level sets of the Hamiltonian when $p\neq0$.

In Figure \ref{figure:Fold}, we showed a cross-section of the manifold $\mathcal{M}_{\beta_1}$ for system \eref{LV3}  at fixed $x_2$ near a fold. In contrast, Figure \ref{figure:foldvaryp} shows a cross-section of the manifold $\mathcal{M}_{\beta_1}$ for system \eref{LVrotate} at the same fixed value of $x_2$. Note that the fold in $\mathcal{M}_{\bq}$, which gives rise to nonuniqueness in solutions of the inverse problem, persists and deepens for negative values of $p$ but disappears when $p$ exceeds some positive value. Disappearance of the fold at positive values of $p$ implies that on this particular segment of $\mathcal{M}_{\beta_1}$ the inverse problem has a unique solution. However, we suspect that for low values of $y_2$ this branch will overlap with the part of $\mathcal{M}_{\beta_1}$ corresponding to the solutions with signature $[+--+]$, which would yield an enlargement of the $\mathcal{R}_{\rm{RG}}$ region.

\begin{figure}
\centering  
\includegraphics[width=0.7\textwidth]{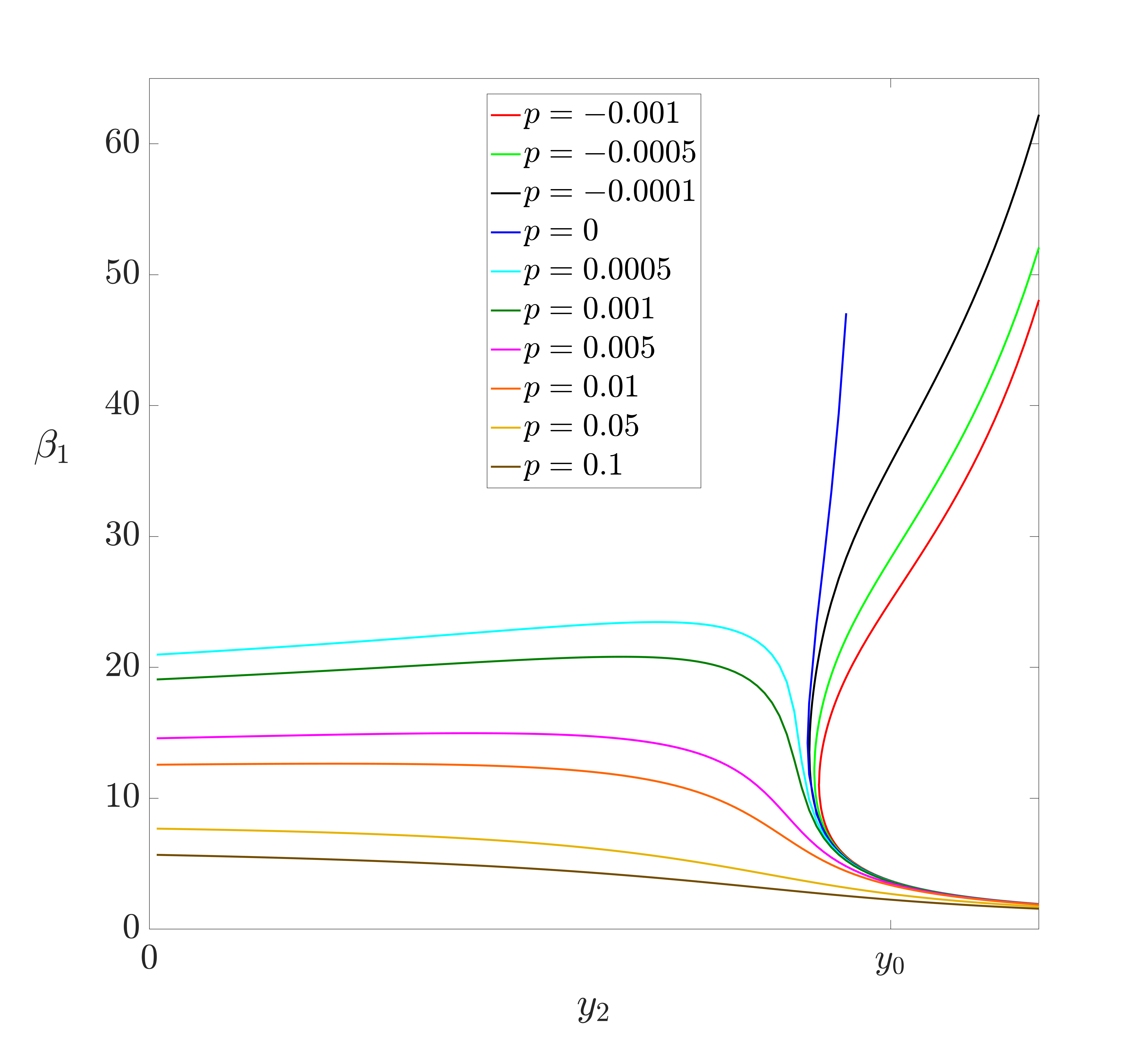}
\caption{Cross-section of the manifold $\mathcal{M}_{\bq}$ at fixed $x_2 = 4.5$ (with $P_0=[1,1]$ and $P_1=[2,1.5]$) for system \eref{LVrotate} with different values of the parameter $p$. Note that folds occur only when $p$ lies below a certain threshold.}\label{figure:foldvaryp}           
\end{figure}

 \subsection{Saturated Lotka-Volterra model}
 
Another natural variation on system \eref{LV3} is the introduction of saturation in the interaction terms in the model.  For example, in a predator-prey model with  predator-prey interaction term $r(x)y$, the rate $r(x)$ can naturally saturate as $x$ increases, due to satiation of the predators.  To implement this modification, we consider the  Lotka-Volterra model
\begin{equation}\label{LV3SaturX}             
\left\{
             \begin{array}{lr}
             \dot{x}=\aq x+\bq xy/(\epsilon x+1), &  \\
             \dot{y}=\bw xy/(\epsilon x+1)+\aw y, &  
             \end{array}
\right.
\end{equation}
where $\epsilon$ is a nonnegative parameter. For $\epsilon=0$ system \eref{LV3SaturX} reduces to system \eref{LV3}, while with $\epsilon>0$ the system includes a saturating dependence of species interaction rate on the value of $x$.

For this system, we first choose data points $P_0=[1,1], P_1=[2,1.5], P_2=[2.45,4]$ and calculate the parameter matrix $A$ for which the trajectory of \eref{LV3SaturX} passes through these data points at times $t=0, t=1, t=2$, respectively. Here $\sigma_A=[+ - + -]$, and from Section \ref{subsection:non-uniqueness:green} we know that  there exist infinitely many solutions with trajectories that overlap when $\epsilon=0$.

We consider what happens to a selection of these solutions as we increase $\epsilon$ from 0.
We choose seven trajectories, three with clockwise (cw) rotation and four with counterclockwise (ccw) rotation. We continue each of these trajectories as a solution to the inverse problem for \eref{LV3SaturX} as $\epsilon$ increases and we visualize the resulting curves of obtained parameter values plotted versus $\epsilon$ in Figure \ref{figure:saturX:2p45:4:eps0to1}. We find that  only two of the seven curves ($\rm{ccw}_0$ and $\rm{cw}_1$) persist over the whole interval  $\epsilon \in [0,1]$.
The other curves each have folds at some $\epsilon \in (0,1)$ , which implies that 
the corresponding solutions are lost before $\epsilon$ reaches 1, in some cases for quite small values of $\epsilon$.

\begin{figure}
\centering  
\includegraphics[width=0.9\textwidth]{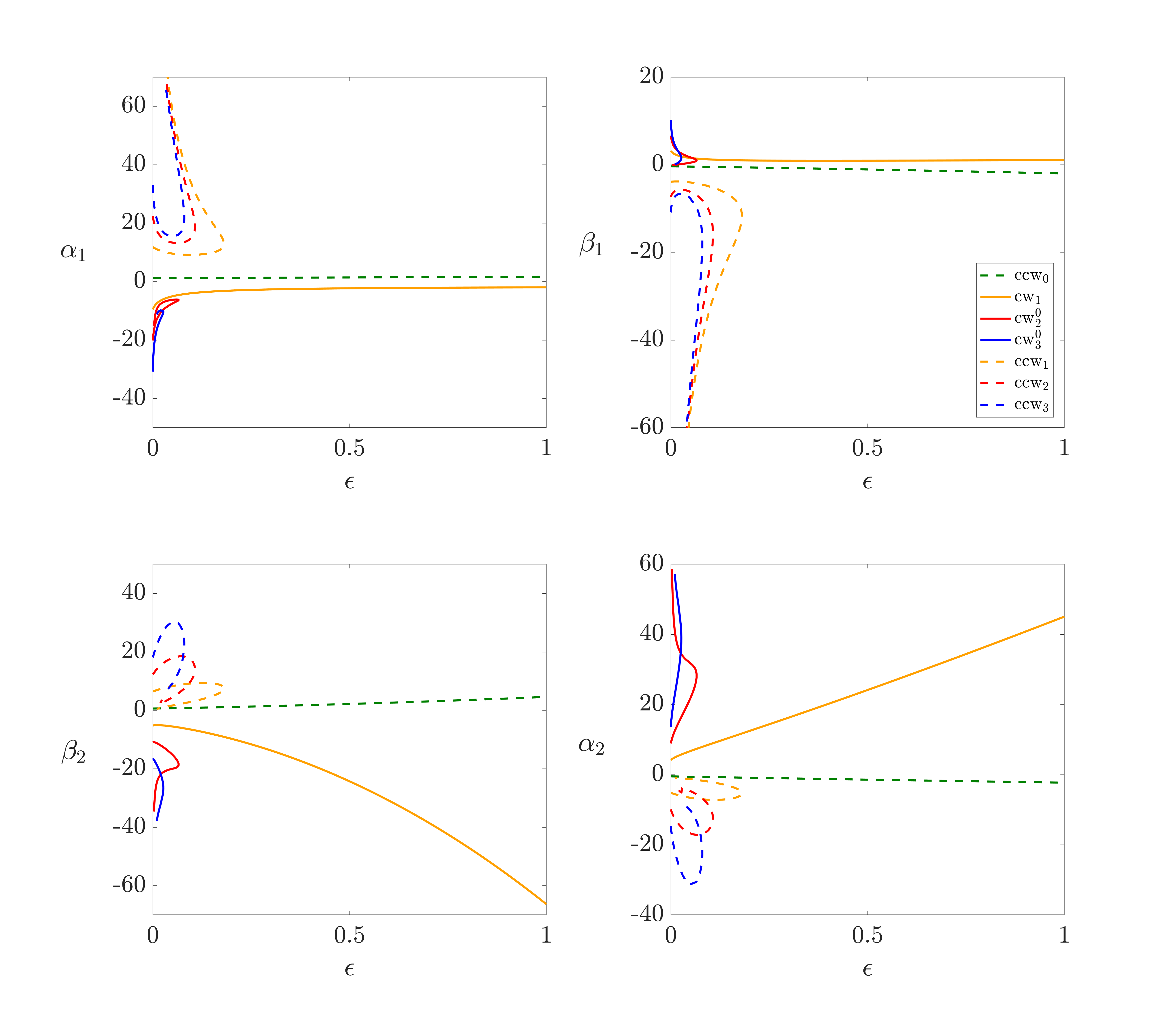}
\caption{Continuation in the direction of increasing $\epsilon$ of seven inverse problem solutions of \eref{LVrotate} corresponding  to data points $P_0=[1,1], P_1=[2,1.5], P_2=[2.45,4]$, starting from $\epsilon=0$. The subscripts in the legend refer to the number of rotations in the trajectory between $t=0$ and $t=1$.}\label{figure:saturX:2p45:4:eps0to1}           
\end{figure}

Interestingly, we do find multiple co-existing inverse problem solutions at $\epsilon = 1$.  The co-existent solutions yield spiral trajectories, some cw and some ccw, some but not all of which exhibit the same number of rotations between passages through the data points.
In Figure \ref{figure:saturX:2p45:4:eps1to0}, we show the continuation of six such solution curves as $\epsilon$ decreases from 1, plotted versus each of the four original model parameters. The dashed green curve and the solid orange curve are exactly the same curves as shown in Figure \ref{figure:saturX:2p45:4:eps0to1}. The other four curves are only defined for $\epsilon\in (0,1]$,  and appear to approach vertical asymptotes,  where one or more parameter values go to $\pm \infty$, at $\epsilon = 0$.
In order to distinguish the $\rm{cw}_2$ and $\rm{cw}_3$ appearing in both Figure \ref{figure:saturX:2p45:4:eps0to1} and Figure \ref{figure:saturX:2p45:4:eps1to0}, we label them with superscripts $0$ or $1$ respectively.


\begin{figure}
\centering  
\includegraphics[width=0.9\textwidth]{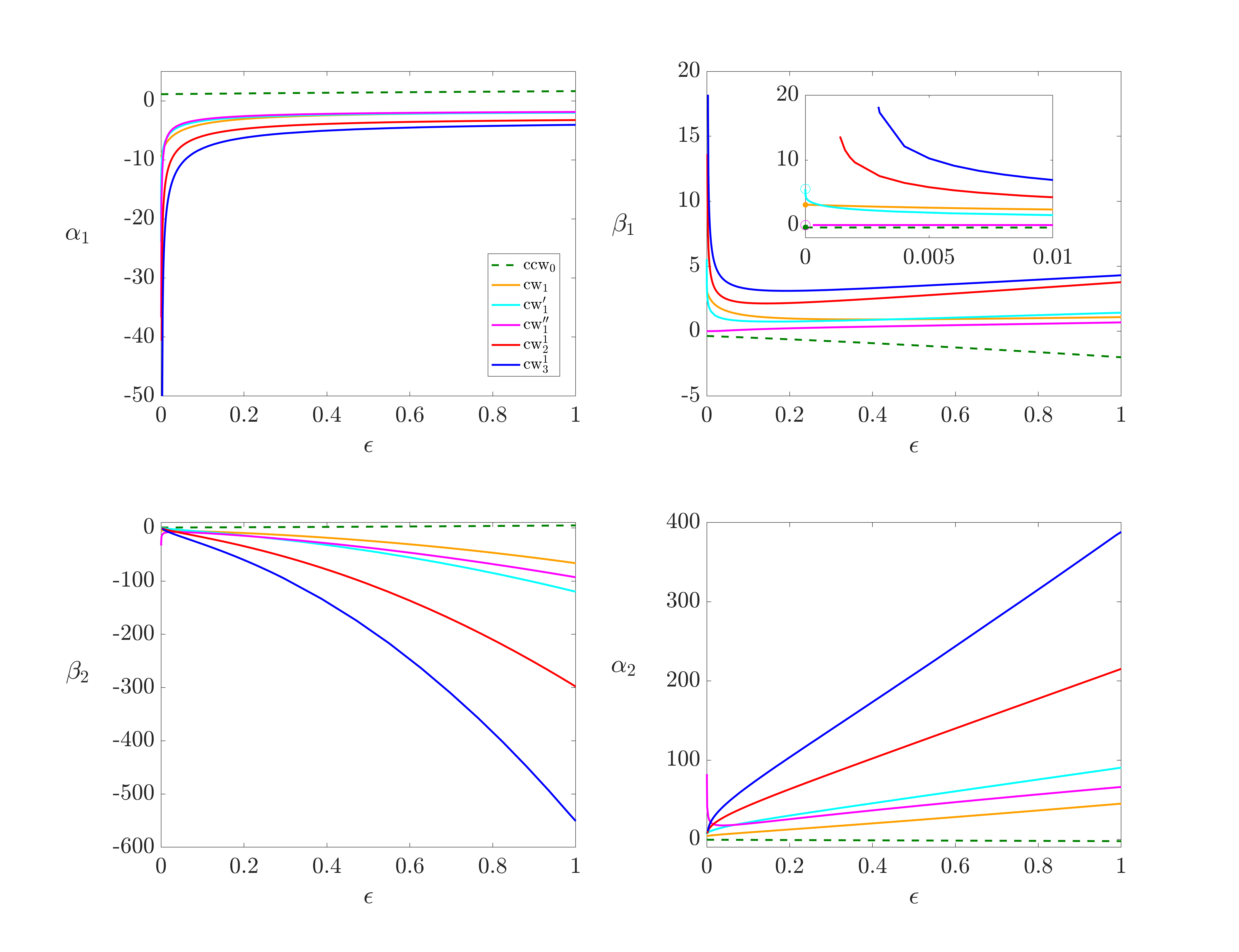}
\caption{Continuation in the direction of decreasing $\epsilon$ of six inverse problem solutions of \eref{LVrotate} corresponding  to data points $P_0=[1,1], P_1=[2,1.5], P_2=[2.45,4]$, starting from $\epsilon=1$. The subscripts in the legend refer to the number of rotations in the trajectory between $t=0$ and $t=1$. }\label{figure:saturX:2p45:4:eps1to0}           
\end{figure}


In Figure \ref{figure:saturX:2p45:4:eps1to0}, there are different alternative solutions that all have one rotation between data points (examples of trajectories corresponding to selected points on the curves in Figure \ref{figure:saturX:2p45:4:eps0to1} and Figure \ref{figure:saturX:2p45:4:eps1to0} with $\epsilon=0.01$, including these one-rotation solutions, are shown in Figures 1 and 2 in the Supplementary material). 
 These alternative solutions are indicated as curves $\rm{cw}_{1}$(orange), $\rm{cw}_{1}^{\prime}$(cyan) and $\rm{cw}_{1}^{\prime\prime}$(magenta). This non-uniqueness of spirals with the same number of rotations between data points does not necessarily appear for all choices of $P_2$. In the Supplementary material, we show the curves generated with the same $P_0$ and $P_1$ but $P_2=[2.45,0.5]$, and show in Supplementary Figure 3 that in this case, there are three curves that persist over the whole interval $\epsilon\in[0,1]$, while we no longer have non-uniqueness of spirals with the same number of rotations. This example implies that the $P_n$-diagram for the saturated Lotka-Volterra system \eref{LV3SaturX} displays a richer variety of non-uniqueness of solutions of the inverse problem than, for example, arises in linear ODE systems.
 

\section{Discussion}

Estimation of model parameters from data is a key step in modeling physical systems.  The usual approaches to this task involve identifying the parameter set that minimizes some measure of deviation of model trajectories from known data points or deriving a probability distribution on parameter space \cite{smith2013uncertainty,tarantola2005inverse}.  The latter approach can be used sequentially to update the probability distribution, or the parameter estimate together with estimates of uncertainty, as new measurements are obtained \cite{lillacci2010}. 
In this work, we take an alternative approach in which, instead of fixing the data points and the associated inverse problem, we consider a whole space of possible data values and the corresponding class of inverse problems.  We show that this approach, first developed for linear \cite{stanhope2017} and affine \cite{duan2020identification} systems, extends easily to linear-in-parameters models, exemplified by the Lotka-Volterra system and associated perturbations.  In this approach, we obtain a manifold $\mathcal{M}$ over the space of data coordinates and unknown parameters, which represents the mapping from the former to the latter.

Our results show that if we project the manifold $\mathcal{M}$ onto the coordinates $x_n,y_n$, of the last data point necessary to to eliminate systemic unidentifiability, we obtain a diagram (i.e., $P_2$ or, more generally, $P_n$-diagram) that encompasses the answers to almost all of the questions posed in Section \ref{section:mainresult}, and we can visualize additional information by plotting $\mathcal{M}_{\bq}$ as shown in Figure \ref{figure:3Db1vsP2_ver1}. Among the \emph{similarities} of this diagram with the corresponding diagrams for linear and affine systems are the following: 
\begin{enumerate}
    \item Domains containing data that yields existence or nonexistence of the solution to the inverse problem are simply connected regions with boundaries that are composed of smooth curves.
    \item Domains containing data for which the inverse problem solution has specified characteristics are also simply connected regions with boundaries composed of smooth curves.
\end{enumerate}

The knowledge that the regions of data points for which the inverse problem solutions share the same properties are simply connected allows a modeler to have confidence in conclusions drawn from fitting the model to the data, because these conclusions are robust to slight changes in the data; in particular, if the data used to solve the inverse problem deviates from the true properties of the system being modeled due to measurement error, that need not yield erroneous conclusions.  Care must be taken, however, to understand the structure of the $P_n$-diagram for a system of interest, because some of its regions may be small. Estimating the sizes of these regions analytically has been done for linear systems \cite{stanhope2017}, but performing similar analysis for nonlinear systems remains a challenge.

The main \emph{differences} between this diagram and the corresponding diagrams for linear and affine systems are the following:
\begin{enumerate}
    \item There are folds in the diagram that give rise to regions of non-uniqueness where two \emph{qualitatively distinct} models fit the same data.
     \item The structure of the diagram is no longer invariant under changes in the location of the points $P_0,P_1$.
     \item Structural stability of the inverse problem solution set may be lost.
\end{enumerate}

In our analysis, the folds referenced in point (i) appeared as fold curves in Figure \ref{figure:P2diagram} and were illustrated in three-dimensional views in Section \ref{sec:vis} above. Such a fold does not occur in the $P_2$-diagrams for linear and affine systems, where the only form of non-uniqueness arises in regions with spiraling solutions and occurs due to countable families of solutions featuring different numbers of rotations in between crossings through the data points. This form of non-uniqueness occurs in the LV system as well and it was described as the trivial nonuniqueness in Section \ref{subsection:non-uniqueness:green}.

Regarding point (ii), in our construction of the $P_2$-diagram we have focused on the case in which $P_0$ and $P_1$ are located so that $x_0<x_1$ and $y_0<y_1$. 
Figure 10 in Section \ref{sec:altP2} shows the diagram corresponding to other choices of the relative location of $P_0$ and $P_1$. Although there are similarities between those diagrams and that of Figure \ref{figure:P2diagram} in the main text, for example in the location of the lines $\mathcal{C}_{\bq}$ and $\mathcal{C}_{\bw}$, the differences are in the appearance of regions for new signature types and in the location of the nonexistence region $\mathcal{R}_{\rm{NE}}$ and the boundary curves $\mathcal{C}_{\aq}$ and $\mathcal{C}_{\aw}$. Even if the inequalities $x_0<x_1$ and $y_0<y_1$ are preserved, however, we have shown that the position of the intersection point of the curves $\mathcal{C}_{\aq}$ and $\mathcal{C}_{\aw}$ is sensitive to the location of $P_0$ and $P_1$, and that the region $\mathcal{R}_{\rm{B}_4}$ in the diagram may not survive when the ratios $x_1/x_0$ and $y_1/y_0$ are too large.

Finally, we briefly illustrated point (iii) in Section \ref{sec:noncon}.  With the introduction of nonlinearity comes new opportunities for changes in solution structure with changes in the vector field.  For example, arbitrarily small perturbations to a conservative system can make the perturbed system non-conservative.  The inverse problem solutions that we considered for the conservative LV system did continue along smooth curves with the introduction of a saturating perturbation, although the interval of persistence was small in some cases.  On the other hand, solutions present for arbitrarily small, positive values of the perturbation parameter did not necessarily persist when the perturbation was completely eliminated.  The loss of these solutions corresponded to parameter values along a solution branch diverging to infinity.  This outcome is analogous to what we observed along some curves in the $P_2$-diagram with certain variations of data point locations.

We made several simplifying assumptions in this study.  One central assumption was that  the number of data points provided was the smallest amount capable of identifying the system, i.e., half of the number of unknown parameters (six), such that the number of given pieces of information about trajectory coordinates was exactly equal to six.  Clearly, the need to fit fewer data points would represent weaker constraints on inverse problem solutions leading to unidentifiability, while the inclusion of more data, while possibly eliminating non-uniqueness, would add constraints or greatly expand the region of non-existence for a strict inverse problem solution in the absence of measurement errors.  
Another central assumption was that the intervals between times when trajectories cross through the data points were equal.  This is a convenient assumption for analytical treatment of linear and affine systems, so we kept it here for consistency.  The qualitative properties observed in our study, however, do not depend on this assumption.  For example, for any time of passage $t>0$ between $P_1$ and $P_2$, there will always be positions in the $P_2$-diagram such that as $P_2$ approaches such a position, the parameter values must go to $\pm \infty$ to maintain compatibility with the time of passage from $P_0$ to $P_1$.  Nonetheless, there remain many open directions related to the ideas in this study, including rigorous analysis of our numerical observations and extension of the ideas in this work to more general classes of nonlinear systems.

\section{Acknowledgements}
JER was partially supported by NSF award DMS 1951095.

\section*{References}
\bibliography{affine}
 
\end{document}


\title[Parameter identification from single trajectory data]{Supplementary Materials: Parameter identification from single trajectory data:  from linear to nonlinear}
\author{X Duan$^1$, J E Rubin$^1$, D Swigon$^1$}
\address{$^1$ Department of Mathematics, University of Pittsburgh, Pittsburgh, PA 15260}
\ead{xid33@pitt.edu, jonrubin@pitt.edu, swigon@pitt.edu}

\section{Numerical methods used in the study}\label{appendix:numerical}
 

In cases when we do not have explicit expressions for the boundary curves in the $P_2$-diagram (Figures 3, 5, and 8 of the main text), we apply numerical methods to generate them. This step can be achieved by solving the boundary value problem (BVP) 
\[
\left\{
             \begin{array}{lr}
             \dot{x}=\aq x+\bq xy, &  \\
             \dot{y}=\bw xy+\aw y, &  \\
             x(i)=x_i, i=0,1,2, &\\
             y(i)=y_i, i=0,1,2, &\\
             \end{array}
\right.
\]
where $P_i=(x_i,y_i),i=0,1,2$ are given, and the unknowns are not the functions $x(t)$ and $y(t)$ themselves but the parameters $\aq, \aw, \bq, \bw$. For fixed $P_0$ 
and $P_1$
and any given  $P_2=(x_2,y_2)$ in the first quadrant, we solve the BVP to obtain a numerical estimation of these parameters for this $P_2$. We then fix $x_2$ and vary $y_2$ (or vice versa), and estimate the corresponding parameters using a continuation tool; we choose to use the program XPP-AUTO \cite{ermentrout2002simulating,doedel1981auto}.  We terminate the numerical procedure when the absolute values of the parameter values (for example, $\bq$ or $\bw$) reach some large level. For example, the fold in Figure 7 
is sketched by fixing $x_2>x_1^2/x_0$ and varying $y_2$. We apply this continuation method for several $x_2$ or $y_2$ values and make a grid in the first quadrant near the region in Figure 3(c), 
from which we sketch the structures shown in Figures 8 and 9. 

For the curves $\mathcal{C}_{\alpha_i}$ and $\mathcal{C}_{\beta_i}, i=1,2,$ we again apply XPP-AUTO to numerically solve the BVP, but we fix the corresponding parameter to be zero and vary $x_2$ or $y_2$ (depending on whether the curve can be parametrized in $x$ or in $y$, respectively).
For the fold curves $\mathcal{C}_{\rm{f}_i}, i=1,2$, we use the two-parameter-continuation tool in XPP-AUTO, beginning with a $P_2$ at the fold position for some parameter set (for example, the blue dot in Figure 7(a)), 
and using $y_2$ as the first parameter and $x_2$ as the second for continuation to generate $\mathcal{C}_{\rm{f}_1}$ (and vice versa for $\mathcal{C}_{\rm{f}_2}$) and estimate the parameter values along these curves. 

With this method, 
we can sketch all of the curves and regions in Figures 8 and 9. 
We then project the surface onto the $(x,y)-$plane and sketch the $P_2$-diagram as in Figures 3 and 5. 
For Figure 6, 
we also use  XPP-AUTO, but with a polar grid centered at $P_0$ instead of Cartesian grid. We start at some $P_2$ in the region $(x_0,x_1)\times (0,y_0)$, and apply the continuation method along the straight line connecting the starting point and $P_0.$ Notice that the surface in Figure 6
over the region $(0,x_0)\times(y_0,y_1)$ is different from that in Figure 8. 
This difference results from the trivial non-uniqueness we mentioned in Subsection 4.1: 
we start at some $P_2$ in the region $(x_0,x_1)\times (0,y_0)$ where  the greatest possible transit time along the periodic orbit, under the constraint that the passage times from $P_0$ to $P_1$ and from $P_1$ to $P_2$ are both $1$, occurs for a set of parameter values with the sign signature $\sigma_A=[- + - +]$. As  we continue the solution curve and pass through $P_0$ into 
the region $(0,x_0)\times(y_0,y_1)$, a switch occurs such that the solution with $\sigma_A=[+ - + -]$ yields the greatest transit time, yet the continuation method does not permit a jump between different $\sigma_A$ solution families.

\section{Periodic orbit geometry}
\label{Appendix:periodic:orbit:geometry}


First, note that when $\beta_1\beta_2\neq 0$, system (4.1) has two equilibrium points, $(0,0)$ and 
\[
P_*=(x_*,y_*)=\left(-\frac{\aw}{\bw}, -\frac{\aq}{\bq}\right),
\]
which does not have to be in the first quadrant.

Second, it is easy to check that system (4.1) is conservative, with the Hamiltonian function
\[
H(x,y, A)=\aw \ln x+\bw x-\aq \ln y-\bq y=\bw(x-x_*\ln x)-\bq(y-y_*\ln y).
\]
Since the level sets of the Hamiltonian represent orbits of system (4.1), one can derive various relations between the constants in the model from the data. In particular, for a trajectory passing through any two points $(P_0,P_1)$ (regardless of the timing or order of the passage) we have $H(x_0,y_0,A)=H(x_1,y_1,A)$, which implies the following relation between the point coordinates and the ratio $r$ of the constants $\bw$ and $\bq$:
\[
 r=\frac{\bw}{\bq}=\frac{y_1-y_0-y_*\ln\left(\frac{y_1}{y_0}\right)}{x_1-x_0-x_*\ln\left(\frac{x_1}{x_0}\right)}.
\]

System (4.1) 
thus supports periodic orbits as level sets of the Hamiltonian, but they have a certain geometry.  Specifically, we have the following result:
\begin{lemma}\label{lemma:convexorbit} 
The region enclosed by any periodic orbit of system (4.1) 
is convex.
\end{lemma}
\begin{proof}
Without loss of generality we suppose that $H(x,y, A)=\bw(x-x_*\ln x) -\bq (y-y_*\ln y)=K$ is a periodic orbit with $\bq>0$ and $\bw<0$, so the flow proceeds clockwise and the region enclosed by this orbit is $\{(x,y) : H(x,y,A)\geq K\}$. Let $(x_a, y_a)$ and $(x_b, y_b)$ be any two points in the region, then for any $t\in[0,1]$,
$$H\left((1-t)\left[ \begin{array}{cc}   x_a \\ y_a \end{array}\right]
+t\left[\begin{array}{cc} x_b \\ y_b\end{array}\right]\right)$$
$$=\bw\left( (1-t)x_a+t x_b-x_*\ln \left((1-t)x_a+tx_b\right)\right) -\bq \left( (1-t)y_a+t y_b-y_*\ln \left((1-t)y_a+ty_b\right)\right)$$
$$\geq\bw\left( (1-t)x_a+t x_b-x_* \left((1-t)\ln(x_a)+t\ln(x_b)\right)\right) -\bq \left( (1-t)y_a+t y_b-y_* \left((1-t)\ln(y_a)+t\ln(y_b)\right)\right)$$
$$=(1-t)H(x_a, y_a, A)+tH(x_b, y_b, A)\geq K.$$
\end{proof}

\section{Domain boundaries}\label{Appdx:domain:bdries}

In this section, we fix the position of the first two data points $P_0$ and $P_1$ such that $x_1>x_0$ and $y_1>y_0$, and explore explicit functions that determine some of the domain boundaries or intersection points of those boundaries.

\subsection{$\mathcal{C}_{\bq}$ or $\mathcal{C}_{\bw}$}
When $\bq=0$, the first equation in system (4.1) 
becomes $\dot{x}=\aq x$, which is independent of $y$, and the whole system can be solved analytically. Its solution is
\[
\left\{
             \begin{array}{rl}
              x(t)=&\rme^{\aq t}x_0,\\
              y(t)=&y_0 \exp\left(\aw t+ \frac{\bw x_0}{\aq}  (\rme^{\aq t}-1) \right).
             \end{array}
\right.
\]
We can also represent the solution curve as a graph over $x$ of the function
\[
\begin{array}{rl}
    y(x)&=y_0\left(\frac{x}{x_0}\right)^{\aw/ \aq}\exp\left(\frac{\bw}{\aq}(x-x_0)\right) \\
     &=y_0\exp\left(\frac{\aw}{\aq}\ln\left(\frac{x}{x_0}\right)+\frac{\bw}{\aq}(x-x_0)\right). 
\end{array}
\]

As $x_2=\rme^{2\aq}x_0=(\rme^{\aq}x_0)^2/x_0=x_1^2/x_0$, the boundary $\mathcal{C}_{\bq}$ is the vertical straight line $x=x_1^2/x_0$ in the $(x,y)-$plane. From the solution formulas above, when $\bq=0$, the parameter matrix can be identified for any point $P_2$ on $\mathcal{C}_{\bq}$ as follows:
\[
A= \left[
\begin{array}{ccc}
\ln \left(\frac{x_1}{x_0}\right) & 0 & 0\\
0 & \frac{x_0}{(x_1-x_0)^2}\ln\left(\frac{x_1}{x_0}\right)\ln\left(\frac{y_0y_2}{y_1^2}\right) & \frac{x_0}{x_0-x_1}\ln\left(\frac{y_2}{y_1}\right)-\frac{x_1}{x_0-x_1}\ln\left(\frac{y_1}{y_0}\right)
\end{array}
\right].
\]

The point of intersection of $\mathcal{C}_{\bq}$ and $\mathcal{C}_{\aw}$ can now be related directly to $P_0$ and $P_1$ as 
\[
P_{\bq,\aw}=\left(x_1\left(\frac{x_1}{x_0}\right), y_1\left(\frac{y_1}{y_0}\right)^\frac{x_1}{x_0}\right).
\]

The case when $\bw=0$ can be studied similarly, and $\mathcal{C}_{\bw}$ is the horizontal line $y=y_1^2/y_0$. For the intersection point of $\mathcal{C}_{\bw}$ and $\mathcal{C}_{\aq}$ we have
\[
P_{\bw,\aq}=\left(x_1\left(\frac{x_1}{x_0}\right)^\frac{y_1}{y_0}, y_1\left(\frac{y_1}{y_0}\right)\right).
\]




\subsection{$\mathcal{C}_{\aq}$ or $\mathcal{C}_{\aw}$}
\label{aq=0}

When  $\aq=0$, system (4.1) 
reduces to
\[
\begin{array}{rl}
             \dot{x}&=\bq xy  \\
             \dot{y}&=\bw xy+\aw y
\end{array}
\]
and every point on the $x$-axis is a critical point of the system. From the Hamiltonian one can deduce that 
\[
    y =\bq^{-1}(  \aw \ln x+\bw x- C )=r(x-x_* \ln x - \bar{C}).
\]
where knowledge of $P_0$ and $P_1$ determines, to within the parameter $x_*$, the constant $r$, given by
\[  
r =\frac{y_1-y_0}{x_1-x_0-x_*\ln\left(\frac{x_1}{x_0}\right)},
\]
and the constant $\bar{C}$, given by 
\[
\bar{C} = 
\frac{x_0 y_1-y_0 x_1}{y_1-y_0} -x_* \frac{(\ln x_0)y_1 -y_0(\ln x_1)}{y_1-y_0}.
\]
Note that if $0 < x_*<(x_1-x_0)/\ln(x_1/x_0)$ then $r>0$ and
$y\to \infty$ as $x \to 0$ or $x \to \infty$. 

The location of the point $x_2$ can then be determined, for any $x_*$, by quadrature as the value for which the following equation holds:
\begin{equation}\label{integrala1=0}
\int_{x_1}^{x_2} \frac{dx}{x(x-x_* \ln x - \bar{C})} = \bw =  \int_{x_0}^{x_1} \frac{dx}{x(x-x_* \ln x - \bar{C})}.  
\end{equation}

Unfortunately, the integrals cannot be evaluated explicitly and hence the corresponding $P_2$ must be computed numerically. The curves $\mathcal{C}_{\alpha_i}, i=1,2$ are shown in Figure 5 
and discussed in Section 4.2.1. 

\subsection{Intersection of $\mathcal{C}_{\aq}$ and $\mathcal{C}_{\aw}$}
When $\aq=\aw=0,$ the solution of the system (4.1) 
satisfies $\bw\dot{x} = \bq\dot{y} $ and hence lies on a straight line passing through $P_0$ and $P_1$. Thus, the first equation of the system (4.1) 
can be rewritten as a Bernoulli differential equation for  $x$:
\[
    \dot{x}- C\bq x=\bw x^2,
\]
where $C=y_0-x_0\frac{y_1-y_0}{x_1-x_0}$ is the $y-$intercept of the trajectory line and $\bw/\bq = \frac{y_1-y_0}{x_1-x_0}$ is its slope. The 
solution to the Bernoulli equation depends on $C, \bq$ and $\bw$, and is given by standard methods as
\[
\frac{1}{x(t)} = \frac{1}{x_0} e^{-C\bq t} + \frac{y_1-y_0}{x_1y_0 - x_0y_1} (e^{-C\bq t}-1)
\]
where we used the identities $\frac{\bw}{C\bq} = \frac{y_1-y_0}{C(x_1-x_0)} = \frac{y_1-y_0}{x_1y_0 - x_0y_1}$.
By substituting 1 and 2 for $t$ we obtain expressions relating $C\bq$ and $x_2$ to the coordinates $x_0,y_0,x_1,y_1$: 
\[
\frac{1}{x_1} = \frac{1}{x_0} e^{-C\bq} + \frac{y_1-y_0}{x_1y_0 - x_0y_1} (e^{-C\bq}-1),
\]
\[
\frac{1}{x_2} = \frac{1}{x_0} e^{-2C\bq} + \frac{y_1-y_0}{x_1y_0 - x_0y_1} (e^{-2C\bq}-1).
\]
The first implies $e^{-C\bq} = (y_1x_0)/(x_1y_0)$ while together they imply $ e^{-C\bq} = \frac{1/x_2-1/x_1}{1/x_1-1/x_0}$.
Therefore, when $\aq=\aw=0$, $P_2$ depends on $P_0$ and $P_1$ as follows:
\[
    P_{\aq,\aw}=\left(\frac{y_0x_1^2}{x_1y_0+x_0y_1-x_1y_1}, \frac{x_0y_1^2}{x_1y_0+x_0y_1-x_1y_1} \right).
\]
 or, equivalently,
\[
 P_{\aq,\aw}=\left(x_1\left(\frac{x_1}{x_0}\right)\left(\frac{x_1}{x_0}+\frac{y_1}{y_0}-\frac{x_1y_1}{x_0y_0}\right)^{-1},
 y_1\left(\frac{y_1}{y_0}\right)\left(\frac{x_1}{x_0}+\frac{y_1}{y_0}-\frac{x_1y_1}{x_0y_0}\right)^{-1} \right).
\]

This result can be generalized to the case when $\aq=\aw=\alpha$ \cite{varma1977exact}. In that case, the trajectory of the system (4.1) 
satisfies $\bw(\dot{x} - \alpha x) = \bq(\dot{y} -\alpha y) $, which can be integrated as
\[
\bw x(t) = \bq y(t) - C e^{\alpha t}
\]
with $C$ a constant of integration. By substituting 0 and 1 for $t$ one finds that $C = \bq y_0 -\bw x_0 $ and $e^{\alpha} = (\bq y_1-\bw x_1)/(\bq y_0 - \bw x_0)$.
Thus, the first equation of the system (4.1) 
can  be rewritten as a Bernoulli differential equation for $x$:
\[
    \dot{x}- C\bq e^{\alpha t} x=\bw x^2.
\]
Standard methods yield the following solution:
\[
\left(\frac{1}{x(t)} - \frac{1}{x_0}\right) \exp\left(\frac{C\bq }{\alpha}e^{\alpha t}\right) = -\int_0^t \bw \exp\left(\frac{C\bq }{\alpha}e^{\alpha s}\right) ds,
\]
\[
\frac{1}{x(t)} = \frac{1}{x_0} -  \int_0^t \bw \exp\left(\frac{C\bq }{\alpha}(e^{\alpha s} - e^{\alpha t}) \right) ds.
\]
A similar expression can be obtained for $y(t)$ as well.

\section{Other figures for saturated Lotka-Volterra model}

\begin{figure}
\centering  
\includegraphics[width=1\textwidth]{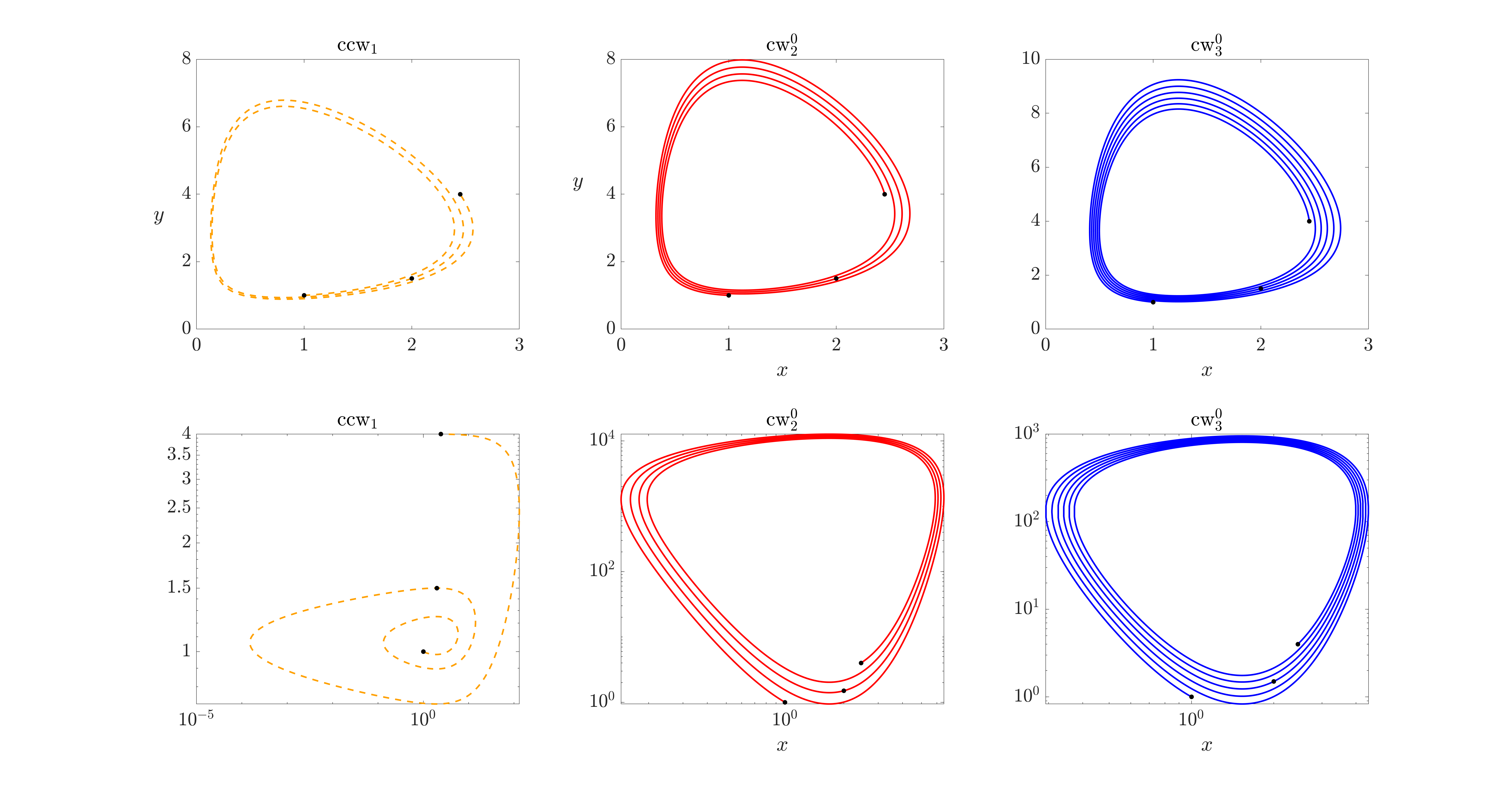}
\caption{$P_0=[1,1], P_1=[2,1.5], P_2=[2.45,4]$, 6 different solutions when $\epsilon=0.01$  shown in three pairs because of the folds. The upper three subplots correspond to the solutions with smaller absolute values of $\beta_1$ while the lower three subplots correspond to the solutions with larger $\beta_1$ values.}\label{figure:saturX:2p45:4:from0:8trajs}          
\end{figure}

\begin{figure}
\centering  
\includegraphics[width=1\textwidth]{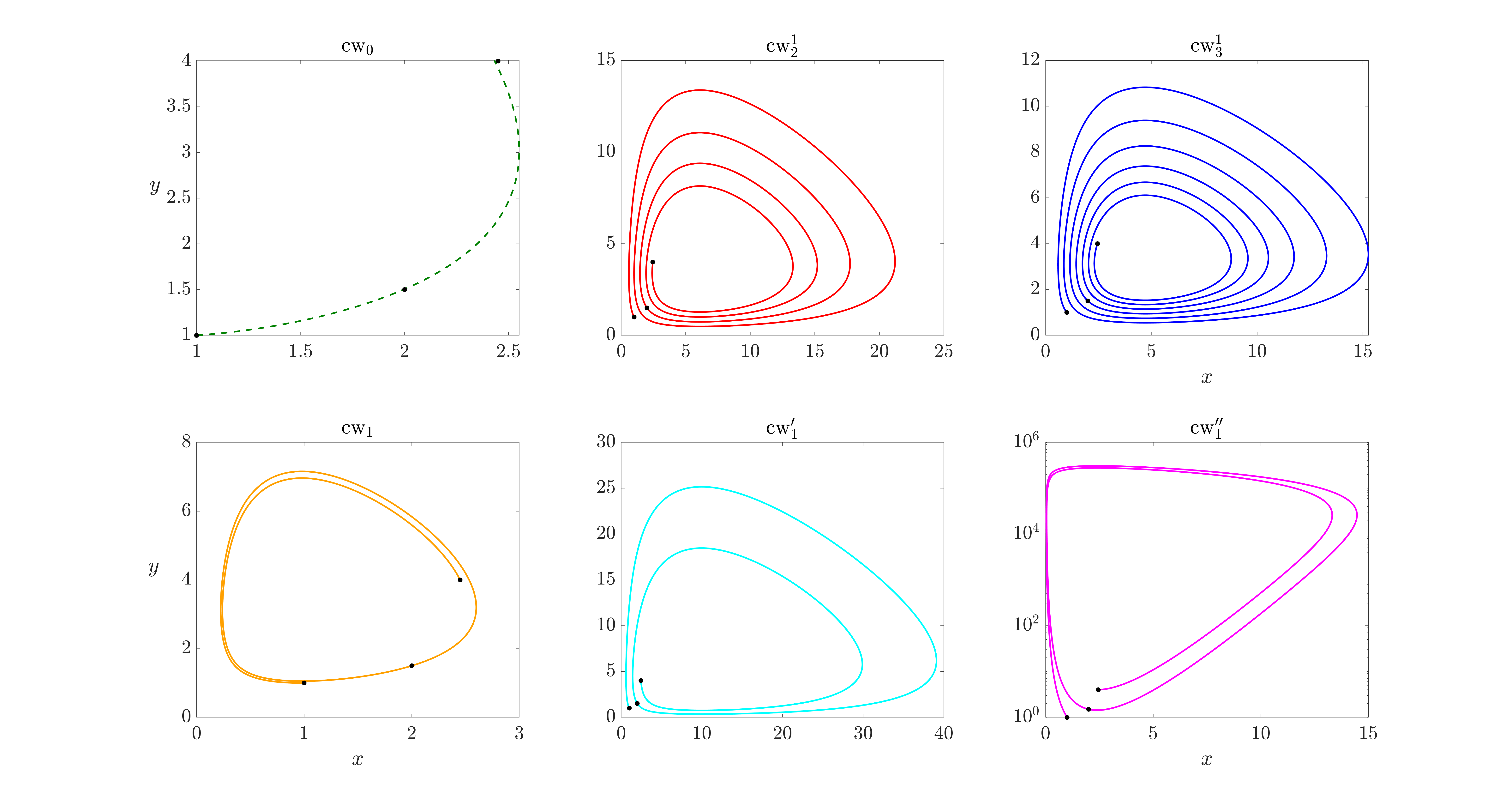} 
\caption{$P_0=[1,1], P_1=[2,1.5], P_2=[2.45,4]$, 6 different solutions when $\epsilon=0.01$.}\label{figure:saturX:2p45:4:from1:6trajs}          
\end{figure}

In Figure 11 and Figure 12 of the paper, we have shown different curves of parameter values versus $\epsilon$. Here we fix $\epsilon=0.01$ and sketch some of the corresponding trajectories in Figure \ref{figure:saturX:2p45:4:from0:8trajs} and Figure \ref{figure:saturX:2p45:4:from1:6trajs}.

\begin{figure}
\centering  
\includegraphics[width=1\textwidth]{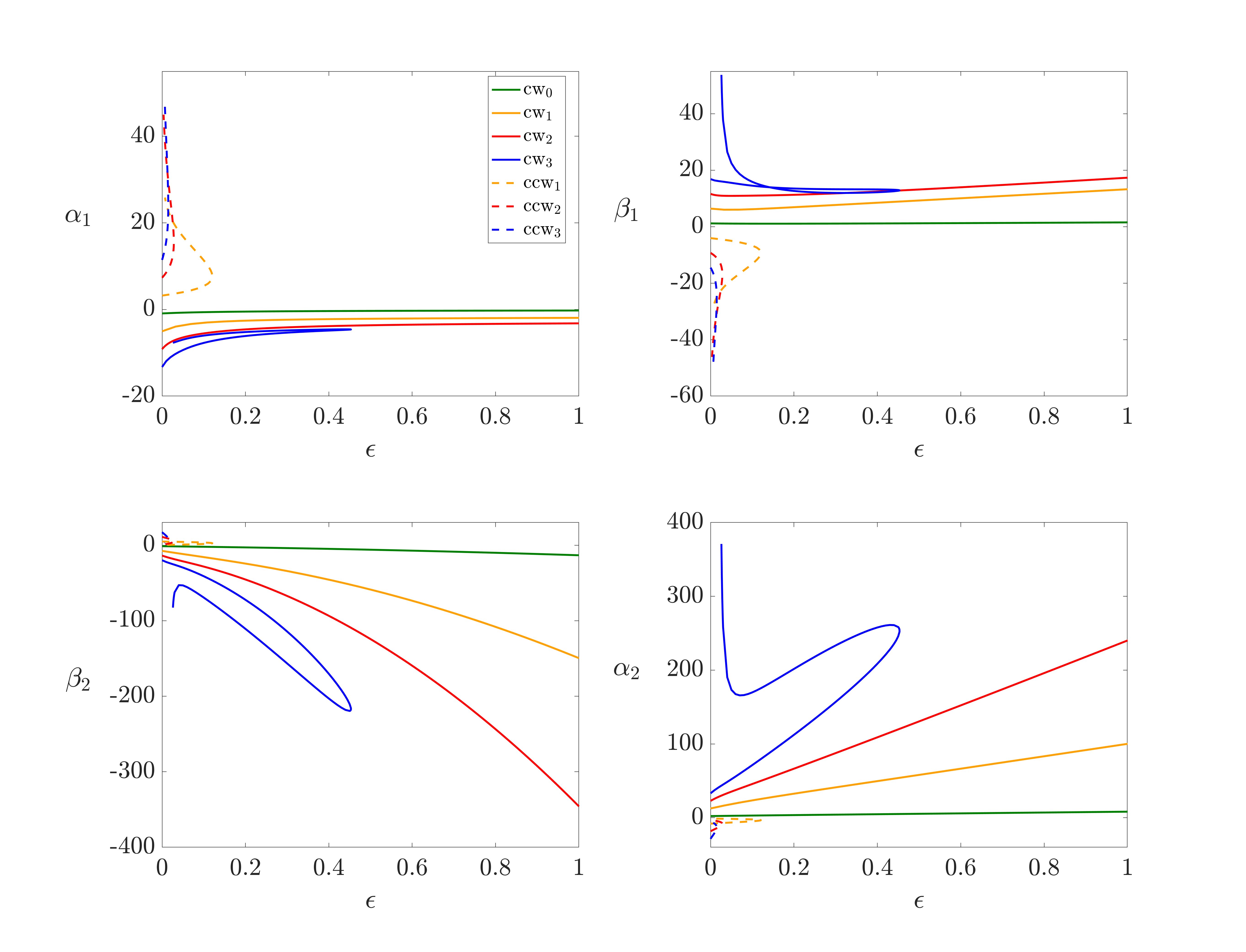}
\caption{$P_0=[1,1], P_1=[2,1.5], P_2=[2.45,0.5]$, starting from $\epsilon=0$.}\label{figure:saturX:2p45:0p5:eps01}           
\end{figure}

To show the curves for other choices of $P_2$, we sketch Figure \ref{figure:saturX:2p45:0p5:eps01} which no longer has several curves all with same 
number of rotations between the data points.

\bigskip

\section*{References}
 
\bibliographystyle{siamplain}
\bibliography{Affine}